\documentclass[11pt]{amsart}
\usepackage{amsmath}
\usepackage{amsxtra}
\usepackage{amscd}
\usepackage{amsthm}
\usepackage{amsfonts}
\usepackage{amssymb}
\usepackage{eucal}

\newtheorem{thm}{Theorem}[section]
\newtheorem{conj}{Conjecture}

\newtheorem{lem}[thm]{Lemma}

\theoremstyle{remark}
\newtheorem{remark}[thm]{Remark}

\theoremstyle{definition}
\newtheorem{definition}[thm]{Definition}

\numberwithin{equation}{section}

\newcommand{\bean}{\begin{eqnarray}}
\newcommand{\eean}{\end{eqnarray}}
\newcommand{\be}{\begin{displaymath}}
\newcommand{\ee}{\end{displaymath}}
\newcommand{\bea}{\begin{eqnarray*}}   
\newcommand{\eea}{\end{eqnarray*}}

\newcommand{\thmref}[1]{Theorem~\ref{#1}}
\newcommand{\secref}[1]{Section~\ref{#1}}
\newcommand{\lemref}[1]{Lemma~\ref{#1}}

\newcommand{\conjref}[1]{Conjecture~\ref{#1}}
\newcommand{\nc}{\newcommand}
\nc{\on}{\operatorname}
\nc{\ch}{\mbox{ch}}
\nc{\Z}{{\mathbb Z}}
\nc{\C}{{\mathbb C}}
\nc{\pone}{{\mathbb C}{\mathbb P}^1}
\nc{\pa}{\partial}
\nc{\F}{{\mathcal F}}
\nc{\arr}{\rightarrow}
\nc{\larr}{\longrightarrow}
\nc{\al}{\alpha}
\nc{\ri}{\rangle}
\nc{\lef}{\langle}
\nc{\W}{{\mathcal W}}
\nc{\la}{\lambda}
\nc{\ep}{\epsilon}
\nc{\su}{\widehat{{\mathfrak s}{\mathfrak l}}_2}
\nc{\sw}{{\mathfrak s}{\mathfrak l}}
\nc{\g}{{\mathfrak g}}
\nc{\h}{{\mathfrak h}}
\nc{\n}{{\mathfrak n}}
\nc{\N}{\widehat{\n}}
\nc{\G}{\widehat{\g}}
\nc{\De}{\Delta}
\nc{\gt}{\widetilde{\g}}
\nc{\Ga}{\Gamma}
\nc{\one}{{\mathbf 1}}
\nc{\z}{{\mathfrak Z}}
\nc{\La}{\Lambda}
\nc{\wt}{\widetilde}
\nc{\wh}{\widehat}
\nc{\cri}{_{\kappa_c}}
\nc{\sun}{\widehat{\sw}_N}
\nc{\si}{\sigma}
\nc{\el}{\ell}
\nc{\bi}{\bibitem}
\nc{\om}{\omega}
\nc{\ol}{\overline}
\nc{\ds}{\displaystyle}
\nc{\dzz}{\frac{dz}{z}}
\nc{\Res}{\on{Res}}
\nc{\mc}{\mathcal}
\nc{\Cal}{\mathcal}
\nc{\bb}{{\mathfrak b}}
\nc{\ot}{\otimes}
\nc{\R}{{\mc R}}
\nc{\yy}{{\mc Y}}
\nc{\ga}{\gamma}

\nc{\us}{\underset}
\nc{\opl}{\oplus}
\nc{\beq}{\begin{equation}}
\nc{\Fq}{{\mathbb F}_q}
\nc{\Mq}{{\mathcal M}}
\nc{\Rep}{\on{Rep}}
\nc{\sssec}{\subsubsection}
\nc{\ssec}{\subsection}
\nc{\lan}{\langle}
\nc{\ran}{\rangle}

\nc{\D}{\mathcal D}
\nc{\Vect}{\on{Vect}}
\nc{\ghat}{\G}
\nc{\T}{\mc T}
\nc{\Tloc}{\T^\g_{\on{loc}}}
\nc{\vac}{|0\ran}
\nc{\Wick}{{\mb :}}
\nc{\mb}{\mathbf}
\nc{\delz}{\partial_z}
\nc{\K}{{\cali K}}
\nc{\cali}{\mathcal}
\nc{\li}{\mathfrak l}
\nc{\lt}{\widetilde{\li}}
\nc{\astar}{a^*}
\nc{\cA}{{\mc A}}
\nc{\ka}{\kappa}

\nc{\OO}{{\mc O}}
\nc{\AutO}{\on{Aut}\OO}
\nc{\DerO}{\on{Der}\OO}
\nc{\DerpO}{\on{Der}_+\OO}
\nc{\Au}{{\mc A}ut}
\nc{\mf}{\mathfrak}
\nc{\V}{{\mc V}}
\nc{\hh}{\wh{\h}}

\nc{\pp}{{\mathfrak p}}
\nc{\mm}{{\mathfrak m}}
\nc{\rr}{{\mathfrak r}}
\nc{\ket}{\rangle}
\nc{\zz}{{\mathfrak z}}
\nc{\gr}{\on{gr}}
\nc{\Spe}{\on{Spec}}
\nc{\rv}{\rho^\vee}
\nc{\can}{\on{can}}
\nc{\CC}{\on{Op}_G(D))}
\nc{\Op}{\on{Op}_G(D)}
\nc{\MOp}{\on{MOp}_G(D)}
\nc{\Db}{{\mathbb D}}
\nc{\ww}{w}

\nc{\oQl}{\ol{{\mathbb Q}}_\ell}
\nc{\oFq}{\ol{{\mathbb F}}_q}
\nc{\Q}{{\mathbb Q}}
\nc{\Ql}{{\mathbb Q}_\ell}
\nc{\bs}{\backslash}
\nc{\E}{{\mc E}}
\nc{\AD}{{\mathbb A}}
\nc{\M}{{\mc M}}
\nc{\Bun}{\on{Bun}}
\nc{\hl}{h^{\leftarrow}}
\nc{\hr}{h^{\rightarrow}}
\nc{\supp}{\on{supp}}
\nc{\He}{\on{H}}
\nc{\Aut}{\on{Aut}}
\nc{\Ll}{{\mc L}}
\nc{\Coh}{{{\mathcal C}oh}}
\nc{\ovc}{\overset{\circ}}
\nc{\Hav}{\on{H}}
\nc{\Mod}{\on{Mod}}
\nc{\kk}{{\mathfrak k}}
\nc{\vf}{\varphi}
\nc{\Gr}{\on{Gr}}
\nc{\gen}{\on{gen}}
\nc{\IC}{\on{IC}}
\nc{\Jac}{\on{Jac}}

\begin{document}

\title[Recent advances in the Langlands Program]{Recent advances in
  the Langlands Program}

\author{Edward Frenkel}\thanks{Partially supported by grants from the
  Packard Foundation and the NSF}

\address{Department of Mathematics, University of California,
  Berkeley, CA 94720, USA}

\date{December 2002, revised in September 2003. Notes for the lecture
at the ``Current Events'' special session of the AMS meeting in
Baltimore, January 17, 2003.}

\maketitle

\section{Introduction}

\subsection{}

The Langlands Program has emerged in the late 60's in the form of a
series of far-reaching conjectures tying together seemingly unrelated
objects in number theory, algebraic geometry, and the theory of
automorphic forms \cite{La}. To motivate it, recall the classical
Kronecker-Weber theorem which describes the maximal abelian extension
$\Q^{\on{ab}}$ of the field $\Q$ of rational numbers (i.e., the
maximal extension of $\Q$ whose Galois group is abelian). This theorem
states that $\Q^{\on{ab}}$ is obtained by adjoining to $\Q$ all roots
of unity; in other words, $\Q^{\on{ab}}$ is the union of all
cyclotomic fields $\Q(\zeta_N)$ obtained by adjoining to $\Q$ a
primitive $N$th root of unity $\zeta_N$.

The Galois group $\on{Gal}(\Q(\zeta_N)/\Q)$ of automorphisms of
$\Q(\zeta_N)$ preserving $\Q \subset \Q(\zeta_N)$ is isomorphic to the
group $(\Z/N\Z)^\times$ of units of the ring $\Z/N\Z$. Indeed, each
element $m \in (\Z/N\Z)^\times$, viewed as an integer relatively prime
to $N$, gives rise to an automorphism of $\Q(\zeta_N)$ which sends
$\zeta_N$ to $\zeta_N^m$. Therefore we obtain that the Galois group
$\on{Gal}(\Q^{\on{ab}}/\Q)$, or, equivalently, the maximal abelian
quotient of $\on{Gal}(\ol{\Q}/\Q)$, where $\Q$ is an algebraic closure
of $\Q$, is isomorphic to the projective limit of the groups
$(\Z/N\Z)^\times$ with respect to the system of surjections
$(\Z/N\Z)^\times \to (\Z/M\Z)^\times$ for $M$ dividing $N$. This
projective limit is nothing but the direct product of the
multiplicative groups of the rings of $p$--adic integers,
$\Z_p^\times$, where $p$ runs over the set of all primes. Thus, we
obtain that
$$
\on{Gal}(\Q^{\on{ab}}/\Q) \simeq \prod_p \Z^\times_p.
$$

The abelian class field theory gives a similar description for the
maximal abelian quotient $\on{Gal}(F^{\on{ab}}/F)$ of the Galois group
$\on{Gal}(\ol{F}/F)$, where $F$ is an arbitrary global field, i.e., a
finite extension of $\Q$ (number field), or the field of rational
functions on a smooth projective curve defined over a finite field
(function field). Namely, $\on{Gal}(F^{\on{ab}}/F)$ is almost
isomorphic to the quotient $F^\times\bs{\mathbb A}_F^\times$, where
${\mathbb A}_F$ is the ring of {\em ad\`eles} of $F$, a subring in the
direct product of all completions of $F$ (see \secref{prelim}). Here
we use the word ``almost'' because we need to take the group of
components of this quotient if $F$ is a number field, or its profinite
completion if $F$ is a function field.

When $F = \Q$ the ring $\AD_{\Q}$ is a subring of the direct product
of the fields ${\mathbb Q}_p$ of $p$--adic numbers and the field
${\mathbb R}$ of real numbers, and the quotient $F^\times\bs{\mathbb
A}_F^\times$ is isomorphic to ${\mathbb R}_{>0} \times \ds \prod_p
\Z^\times_p$. Hence the group of its components is $\ds \prod_p
\Z^\times_p$, in agreement with the Kronecker-Weber theorem.

One can obtain complete information about the maximal abelian quotient
of a group by considering its one-dimensional representations. The
above statement of the abelian class field theory may then be
reformulated as saying that one-dimensional representations of
$\on{Gal}(\ol{F}/F)$ are essentially in bijection with one-dimensional
representations of the abelian group ${\mathbb A}_F^\times =
GL_1({\mathbb A}_F)$ which occur in the space of functions on
$F^\times\bs {\mathbb A}_F^\times = GL_1(F)\bs GL_1({\mathbb A}_F)$. A
marvelous insight of Robert Langlands was to conjecture that there
exists a similar description of $n$--{\em dimensional representations}
of $\on{Gal}(\ol{F}/F)$. Namely, he proposed that those may be related
to irreducible representations of the group $GL_n({\mathbb A}_F)$
which are {\em automorphic}, that is those occurring in the space of
functions on the quotient $GL_n(F)\bs GL_n({\mathbb A}_F)$. This
relation is now called the {\em Langlands correspondence}.

\subsection{}

At this point one might ask a legitimate question: why is it important
to know what the $n$--dimensional representations of the Galois group
look like, and why is it useful to relate them to things like
automorphic representations?  There are indeed many reasons for
that. First of all, it should be remarked that according to the
Tannakian phylosophy, one can reconstruct a group from the category of
its finite-dimensional representations, equipped with the structure of
the tensor product. Therefore looking at $n$--dimensional
representations of the Galois group is a natural step towards
understanding its structure. But even more importantly, one finds many
interesting representations of Galois groups in ``nature''. For
example, the group $\on{Gal}(\ol\Q/\Q)$ will act on the geometric
invariants (such as the \'etale cohomologies) of an algebraic variety
defined over $\Q$. Thus, if we take an elliptic curve $E$ over $\Q$,
then we will obtain a two-dimensional Galois representation on its
first \'etale cohomology. This representation contains a lot of
important information about the curve $E$, such as the number of
points of $E$ over $\Z/p\Z$ for various primes $p$.

The point is that the Langlands correspondence is supposed to relate
$n$--dimensional Galois representations to automorphic representations
of $GL_n({\mathbb A}_F)$ in such a way that the data on the Galois
side, such as the number of points of $E(\Z/p\Z)$, are translated into
something more tractable on the automorphic side, such as the
coefficients in the $q$--expansion of the modular forms that
encapsulate automorphic representations of $GL_2(\AD_{\Q})$.

More precisely, one asks that under the Langlands correspondence
certain natural invariants attached to the Galois representations and
to the automorphic representations be matched. These invariants are
the {\em Frobenius conjugacy classes} on the Galois side and the {\em
Hecke eigenvalues} on the automorphic side. Let us review them
briefly. To an $n$--dimensional representation $\sigma$ of
$\on{Gal}(\ol{F}/F)$ that is almost everywhere unramified one attaches
a collection of conjugacy classes in the group $GL_n$ for all but
finitely many points $x$ of the curve $X$ if $F = \Fq(X)$ (or primes
of the ring of integers of $F$ if $F$ is a number field). These are
the images of the Frobenius elements in $\on{Gal}(\ol{F}/F)$ (see
\secref{galois}). The eigenvalues of these conjugacy classes are then
represented by unordered $n$--tuples of numbers
$(z_1(\sigma_x),\ldots,z_n(\sigma_x))$. On the other hand, to an
automorphic representation $\pi$ of $GL_n(\AD)$ one attaches a
collection of eigenvalues of the so-called Hecke operators (see
\secref{aut repr}), which may also be encoded by unordered collections
of $n$--tuples of numbers $(z_1(\pi_x),\ldots,z_n(\pi_x))$ for all but
finitely many $x \in X$. If under the Langlands correspondence we have
$$
\pi \; \longleftrightarrow \; \sigma,
$$
then we should have
\begin{equation}    \label{matching}
(z_1(\pi_x),\ldots,z_n(\pi_x)) = (z_1(\sigma_x),\ldots,z_n(\sigma_x))
\end{equation}
for all but finitely many $x \in X$ (up to permutation).

The Frobenius eigenvalues and the Hecke eigenvalues may be converted
into analytic objects called the $L$--functions. Then the matching of
these eigenvalues becomes the statement that the $L$--functions
corresponding to $\pi$ and $\sigma$ are equal.

It is instructive to consider what this matching means in the simplest
example when $F=\Q$ and $n=1$. In this case the Langlands
correspondence comes from the isomorphism between the abelian quotient
of the Galois group $\on{Gal}(\Q^{\on{ab}}/\Q)$ and the group of
components of $\Q^\times\bs{\mathbb A}_{\Q}^\times$. The matching
condition \eqref{matching} then means that the Frobenius element
$\on{Fr}_p$ corresponding to the prime $p \in \Z$ in
$\on{Gal}(\Q^{\on{ab}}/\Q)$ (more precisely, of its quotient
unramified at $p$) goes under this isomorphism to the image of the
ad\`ele $(1,\ldots,1,p,1,\ldots)$, with $p$ being in the factor
$\Q_p^\times$, in the group of components of $\Q^\times\bs{\mathbb
A}_{\Q}^\times$.

Let us represent the latter group as $\underset{\longleftarrow}\lim \;
(\Z/N\Z)^\times$, where $(\Z/N\Z)^\times$ is considered as the Galois
group $\on{Gal}(\Q(\zeta_N)/\Q)$ of the cyclotomic field
$\Q(\zeta_N)$. Then the above statement translates into the statement
as to what is the image of the Frobenius conjugacy class $\on{Fr}_p$
in $\on{Gal}(\Q(\zeta_N)/\Q) = (\Z/N\Z)^\times$. It says that this
image is equal to $p \on{mod} N$; note that this makes sense only if
$p$ does not divide $N$, which is precisely the condition that $p$ is
unramified in $\Q(\zeta_N)$. Since the Frobenius element $\on{Fr}_p$
comes from the automorphism $y \mapsto y^p$ of the finite field
${\mathbb F}_p$, it is not surprising that it corresponds to the
automorphism of $\Q(\zeta_N)$ sending $\zeta_N$ to $\zeta_N^p$.

Thus, we find that the order of $\on{Fr}_p$ in
$\on{Gal}(\Q(\zeta_N)/\Q)$ is equal to the order of $p$ in
$(\Z/N\Z)^\times$. Therefore we may figure out how $p$ splits in the
ring $\Z[\zeta_N]$. If the order of $\on{Fr}_p$ is equal to $d$ and
$(p) = {\mc P}_1 \ldots {\mc P}_r$, where the ${\mc P}_i$'s are prime
ideals in $\Z[\zeta_N]$, then the residue field of each of the ${\mc
P}_i$'s should be an extension of ${\mathbb F}_p$ of degree $d$, and
therefore $r=\varphi(N)/d$, where $\varphi(N) = |(\Z/N\Z)^\times|$. So
we obtain an application of the matching \eqref{matching} in the case
$n=1$ to the problem of the splitting of primes. For instance, if
$N=4$ we obtain that $p$ splits in $\Z[i]$, i.e., $p$ may be
represented as a sum of two squares of integers,
$$p=(a+bi)(a-bi)=a^2+b^2,$$ if and only if $p \equiv 1 \on{mod} 4$,
which is the statement of one of Fermat's theorems (see \cite{Gelbart}
for more details).

Likewise, in the case of two-dimensional Galois representations
$\sigma$ arising from the first \'etale cohomology of an elliptic
curve $E$ over $\Z$, the eigenvalues of $\on{Fr}_p$ (which are
well-defined for all $p$'s that do not divide the conductor of $E$)
contain information about the number of points of $E$ over
$\Z/p\Z$. Suppose now that $\pi$ is an automorphic representation of
$GL_2(\AD)$ attached to $\sigma$ via the Langlands correspondence. One
assigns to $\pi$ in a natural way a modular form on the upper-half
plane (see, e.g., \cite{Gelbart,Murty}). Then the matching condition
\eqref{matching} relates the numbers of points of $E(\Z/p\Z)$ to the
coefficients in the $q$--expansion of this modular form. The existence
of $\pi$ (and hence of the corresponding modular form) now becomes the
statement of the Taniyama-Shimura (or Taniyama-Shimura-Weil)
conjecture that has recenty been proved by A. Wiles and others (it
implies the Fermat's last theorem). One obtains similar statements by
analyzing from the point of view of the Langlands philosophy the
Galois representations coming from other algebraic varieties, or more
general motives.

\subsection{}

While in the number field case the Langlands correspondence has been
established only in special cases such as the one expressed by the
Taniyama-Shimura conjecture, in the function field case the Langlands
conjecture is now a theorem. It has been proved in the 80's by
V. Drinfeld in the case when $n=2$ and recently by L. Lafforgue for an
arbitrary $n$ in a monumental effort for which both of them have been
awarded the Fields Medals.

In these notes we will focus on the Langlands correspondence in the
function field case. We will start by giving a precise formulation of
the Langlands conjecture in this case (see \secref{prelim}). Then in
\secref{proof} we will describe some of the ideas involved in the
proof of this conjecture given by Drinfeld and Lafforgue.

After that we will discuss in \secref{geometric} the {\em geometric
Langlands correspondence}. This is a geometric version of the
Langlands conjecture (available only in the function field case). It
comes from the observation that $n$--dimensional representations of
$\on{Gal}(\ol{F}/F)$, where $F = \Fq(X)$, may be viewed geometrically
as rank $n$ local systems on $X$. Such objects make sense both when
$X$ is a curve over a finite field and over $\C$. Thus, we may
transport the objects on the Galois side of the Langlands
correspondence to the realm of complex algebraic geometry.

It turns out that one may replace the automorphic repsentations of
$GL_n(\AD)$ by geometric objects as well. Those are the {\em Hecke
eigensheaves}, which are certain sheaves on the moduli space of rank
$n$ bundles on the curve. The geometric Langlands conjecture states,
roughly, that to any rank $n$ local system $E$ on $X$ one may
associate a Hecke eigensheaf whose ``eigenvalues'' are expressed in
terms of $E$. The advantage of this formulation is that it makes sense
for curves defined not only over finite fields, but also over the
field of complex numbers. Hence one can use the machinery of complex
algebraic geometry to gain new insights into the nature of the
Langlands correspondence. For example, for $n=1$ the geometric
Langlands conjecture states that to each rank one local system $E$ on
$X$ one can associate a Hecke eigensheaf on the Jacobean variety of
$X$. It turns out that when $X$ is a curve over $\C$ this
correspondence is best understood in the framework of the
Fourier-Mukai transform (see \secref{fm}).

The geometric Langlands conjecture for $GL_n$ has recently been proved
in \cite{FGV,Ga} following earlier works of P. Deligne, V. Drinfeld
\cite{Dr1} and G. Laumon \cite{Laumon:cor1,Laumon:cor2}. In
\secref{geometric} we will give an overview of this proof.

Then in \secref{reductive} we will discuss how to formulate the
Langlands correspondence for an arbitrary reductive group $G$ instead
of $GL_n$, both in the classical and the geometric settings. This is
where one of the most beautiful and mysterious concepts, that of the
{\em Langlands dual group}, enters the subject. Finally, we will
describe in \secref{over C} the work \cite{BD} of A. Beilinson and
V. Drinfeld in which part of this conjecture has been proved for an
arbitrary semisimple Lie group $G$ over $\C$. Their proof uses the
geometry of moduli spaces of bundles on curves as well as
representation theory of affine Kac-Moody algebras.

In the last thirty years the subject of the Langlands correspondence
has evolved into a vast and multifaceted field. In a short review it
is impossible to even glance over its main ideas and trends. In
particular, we will not mention such remarkable recent achievements as
the proof of the local Langlands conjectures for $GL_n$, both in the
function field and number field settings, given in \cite{LST} and
\cite{HT,H}, respectively (in that regard we refer the reader to
Carayol's talk at S\'eminaire Bourbaki \cite{C}).

A number of wonderful reviews of the Langlands Program are available
at present, and I would like to mention some of them. The papers
\cite{Arthur,Arthur2,Gelbart,Knapp,Murty} contain general overviews of
the Langlands Program which are informative and fun to read. The
reports \cite{K} and \cite{HK} give a clear and concise review of the
work of Drinfeld in the case of $GL_2$. For excellent expositions of
Lafforgue's work and the geometric Langlands correspondence we refer
the reader to Laumon's talks at S\'eminaire Bourbaki, \cite{La1} and
\cite{La2}, respectively. Finally, the recent article \cite{Lan}
offers a unique perspective on the field by its Creator ``thirty years
later''.

\section{The statement of the Langlands conjecture in the function
  field case}    \label{prelim}

Let $X$ be a smooth projective curve over a finite field $k = \Fq$
with $q$ elements. We will denote by $\ol{X}$ the corresponding curve
over the algebraic closure $\oFq$ of $\Fq$, $\ol{X} = X
\underset{\on{Spec} \Fq}\times \on{Spec} \oFq$. We will assume that
$\ol{X}$ is connected. Let $F = \Fq(X)$ be the field of rational
functions on $X$, and fix its separable closure $\ol{F}$. We denote by
$G_F$ the Galois group of $F$, i.e., the group of automorphisms of
$\ol{F}$ preserving $F$.

For any closed point $x$ of $X$, we denote by $F_x$ the completion of
$F$ at $x$ and by $\OO_x$ its ring of integers. If we pick a rational
function $t_x$ on $X$ which vanishes at $x$ to order one, then we
obtain isomorphisms $F_x \simeq k_x((t_x))$ and $\OO_x \simeq
k_x[[t_x]]$, where $k_x$ is the residue field of $x$ (the quotient of
the local ring $\OO_x$ by its maximal ideal); this field is a finite
extension of the base field $k$ and hence is isomorphic to ${\mathbb
F}_{q_x}$, where $q_x = q^{\deg x}$. The ring $\AD=\AD_F$ of ad\`eles
of $F$ is by definition the {\em restricted} product of the fields
$F_x$, where $x$ runs over the set $|X|$ of all closed points of
$X$. The word ``restricted'' means that we consider only the
collections $(f_x)_{x \in |X|}$ of elements of $F_x$ in which $f_x \in
\OO_x$ for all but finitely many $x$. The ring $\AD$ contains the
field $F$, which is embedded into $\AD$ diagonally, by taking the
expansions of rational functions on $X$ at all points.

Roughly speaking, the Langlands conjecture states that there is a
bijection between the set of equivalence classes of $n$--dimensional
representations of $G_F$ and the set of equivalence classes of
irreducible automorphic representations of $GL_n(\AD)$, i.e., those
which may be realized in a certain space of functions on the quotient
$GL_n(F)\bs GL_n(\AD)$. In order to make a precise formulation of the
conjecture, we need to explain what all of this means. We will also
restrict ourselves to {\em irreducible} $n$--dimensional
representations of $G_F$. Those will correspond to the so-called {\em
cuspidal} representations of $GL_n(\AD)$. In addition, using the
abelian class field theory, we can without loss of generality consider
only those representations of $G_F$, whose determinant has finite
order, and those representations of $GL_n(\AD)$ whose central
character has finite order (a general representation may be obtained
by tensoring a representation of this type with a one-dimensional
representation). Let us now give the precise definitions.

\subsection{Galois representations}    \label{galois}

Consider the Galois group $G_F$. It is instructive to think of it as a
kind of fundamental group of $X$. Indeed, if $Y \to X$ is a covering
of $X$, then the field $k(Y)$ of rational functions on $Y$ is an
extension of the field $F = k(X)$ of rational functions on $X$, and
the Galois group $\on{Gal}(k(Y)/k(X))$ may be viewed as the group of
``deck transformations'' of the cover. If our cover is unramified,
then this group may be identified with a quotient of the fundamental
group of $X$. Otherwise, this group is isomorphic to a quotient of the
fundamental group of $X$ without the ramification points. The Galois
group $G_F$ itself may be viewed as the group of ``deck
transformations'' of the maximal (ramified) cover of $X$.

Let us pick a point $\ol{x}$ of this cover lying over a fixed point $x
\in |X|$. The subgroup of $G_F$ preserving $\ol{x}$ is called the
decomposition group. If we make a different choice of $\ol{x}$, it
gets conjugated in $G_F$. Therefore we obtain a subgroup of $G_F$
defined up to conjugation. We denote it by $D_x$. This group is in
fact isomorphic to the Galois group $\on{Gal}(\ol{F}_x/F_x)$, and we
have a natural homomorphism $D_x \to \on{Gal}(\ol{k}_x/k_x)$, whose
kernel is called the inertia subgroup and is denoted by $I_x$. The
Galois group $\on{Gal}(\ol{k}_x/k_x)$ has a very simple description:
it contains the (geometric) {\em Frobenius element} $\on{Fr}_x$, which
is defined to be the inverse of the automorphism $y \mapsto y^{q_x}$
of $\ol{k}_x = \ol{\mathbb F}_{q_x}$, and $\on{Gal}(\ol{k}_x/k_x)$ is
equal to the profinite completion of the group $\Z$ generated by this
element.

A homomorphism $\sigma$ from $G_F$ to another group $H$ is called {\em
unramified} at $x$, if $I_x$ lies in the kernel of $\sigma$ (this
condition is independent of the choice of $\ol{x}$). In this case
$\on{Fr}_x$ gives rise to a well-defined conjugacy class in $H$,
denoted by $\sigma(\on{Fr}_x)$.

The group $G_F$ is a profinite group, equipped with the Krull topology
in which the base of open neighborhoods of the identity is formed by
normal subgroups of finite index. Therefore a continuous
finite-dimensional complex representation of $G_F$ necessarily factors
through a finite quotient of $G_F$. To obtain a larger class of Galois
representations we replace $\C$ with the field $\Ql$ of $\ell$--adic
numbers, where $\ell$ is a prime that does not divide $q$. Denote by
$\oQl$ the algebraic closure of $\Ql$. By an $n$--dimensional
$\ell$--{\em adic representation} of $G_F$ we will understand a
homomorphism $\sigma: G_F \to GL_n(\oQl)$ which satisfies the
following conditions:

\begin{enumerate}
\item there exists a finite extension $E \subset \oQl$ of $\Ql$ such
  that $\sigma$ factors through a homomorphism $G_F \to GL_n(E)$,
  which is continuous with respect to the Krull topology on $G_F$ and
  the $\ell$--adic topology on $GL_n(E)$;

\item it is unramified at all but finitely many $x \in |X|$.
\end{enumerate}

Let ${\mc G}_n$ be the set of equivalence classes of irreducible
$n$--dimensional $\ell$--adic representations of $G_F$ such that the
image of $\on{det}(\sigma)$ is a finite group.

Given such a representation, we consider the collection of the
Frobenius conjugacy classes $\{ \sigma(\on{Fr}_x) \}$ in $GL_n(\oQl)$
and the collection of their eigenvalues (defined up to permutation),
which we denote by $\{ (z_1(\sigma_x),\ldots,z_n(\sigma_x)) \}$, for
all $x \in |X|$ where $\sigma$ is unramified. Chebotarev's density
theorem implies the following remarkable result: if two $\ell$--adic
representations are such that their collections of the Frobenius
conjugacy classes coincide for all but finitely many points $x \in
|X|$, then these representations are equivalent.

\subsection{Automorphic representations}    \label{aut repr}

Consider now the group $GL_n(\AD)$. It carries a Haar measure
normalized in such a way that the volume of the subgroup $GL_n(\OO_x)$
is equal to $1$ for all $x \in |X|$. Note that $GL_n(F)$ is naturally
a subgroup of $GL_n(\AD)$. A function $\varphi: GL_n(\AD) \to \oQl$ is
called cuspidal automorphic if it satisfies the following conditions:

\begin{enumerate}
\item it is left $GL_n(F)$--invariant;

\item it is right invariant with respect to an open subgroup of
  $GL_n(\AD)$;

\item there exists an element $a \in \AD^\times$ of non-zero degree
such that $\varphi(ag) = \varphi(g)$ for all $g \in GL_n(\AD)$;

\item let $N_{n_1,n_2}$ be the unipotent radical of the standard
  parabolic subgroup $P_{n_1,n_2}$ of $GL_n$ corresponding to the
  partition $n=n_1+n_2$ with $n_1, n_2>0$. Then
$$
\underset{N_{n_1,n_2}(F)\bs N_{n_1,n_2}(\AD)}\int \varphi(ug) du = 0,
\qquad \forall g \in GL_n(\AD).
$$
\end{enumerate}

Denote the vector space of cuspidal automorphic functions on
$GL_n(\AD)$ by ${\mc C}_{\on{cusp}}$. The group $GL_n(\AD)$ acts on
${\mc C}_{\on{cusp}}$, and under this action ${\mc C}_{\on{cusp}}$
decomposes into a direct sum of irreducible representations. These
representations are called {\em irreducible cuspidal automorphic
representations} of $GL_n(\AD)$. A theorem due to Piatetski-Shapiro
\cite{PS} and Shalika \cite{Sh} says that the each of them enters
${\mc C}_{\on{cusp}}$ with multiplicity one. We denote the set of
equivalence classes of these representations by ${\mc A}_n$.

\medskip

\begin{remark}
If $\pi_1$ and $\pi_2$ are irreducible representations
of $GL_{n_1}(\AD)$ and $GL_{n_2}(\AD)$, respectively, where $n_1+n_2
= n$, then we may extend trivially the representation $\pi_1 \otimes
\pi_2$ of $GL_{n_1} \times GL_{n_2}$ to the parabolic subgroup
$P_{n_1,n_2}(\AD)$ and consider the induced representation of
$GL_n(\AD)$. It is easy to see that irreducible quotients of this
representation cannot be realized in the space ${\mc C}_{\on{cusp}}$
because of the cuspidality condition (4). In fact, this condition is
imposed precisely so as to avoid irreducible quotients of induced
representations.

We also remark that condition (3) is equivalent to the central
character of $\pi$ being of finite order.\qed
\end{remark}

\medskip

Let $\pi$ be an irreducible cuspidal automorphic representation of
$GL_n(\AD)$. One can show that it decomposes into a tensor product 
$$
\pi = \bigotimes_{x \in |X|}{}' \; \pi_x,
$$
where each $\pi_x$ is an irreducible representation of
$GL_n(F_x)$. Furthermore, there is a finite subset $S$ of $|X|$ such
that each $\pi_x$ with $x \in |X|-S$ is {\em unramified}, i.e.,
contains a non-zero vector $v_x$ stable under the maximal compact
subgroup $GL_n(\OO_x)$ of $GL_n(F_x)$. This vector is unique up to a
scalar and we will fix it once and for all. The space $\bigotimes'_{x
\in |X|} \pi_x$ is by definition the span of all vectors of the form
$\bigotimes_{x \in |X|} w_x$, where $w_x \in \pi_x$ and $w_x = v_x$
for all but finitely many $x \in |X|-S$. Therefore the action of
$GL_n(\AD)$ on $\pi$ is well-defined.

Let ${\mc H}_x$ be the space of compactly supported functions on
$GL_n(F_x)$ which are bi-invariant with respect to the subgroup
$GL_n(\OO_x)$. This is an algebra with respect to the convolution
product, which is called the {\em spherical Hecke algebra}. By the
Satake correspondence (see \thmref{satake} below), we have an
isomorphism
\begin{equation}    \label{Hx}
{\mc H}_x \simeq \oQl[z_1^{\pm 1},\ldots,z_n^{\pm 1}]^{S_n}.
\end{equation}
The Hecke algebra ${\mc H}_x$ naturally acts on any irreducible
unramified representation $\pi_x$ of $GL_n(F_x)$ and preserves the
one-dimensional subspace of $GL_n(\OO_x)$--invariant vectors spanned
by $v_x$. Hence ${\mc H}_x$ acts on it via a character, which is
nothing but a collection of non-zero numbers
$(z_1(\pi_x),\ldots,z_n(\pi_x))$ defined up to permutation. We will
call them the {\em Hecke eigenvalues} of $\pi$ at $x$. Thus, to each
irreducible cuspidal automorphic representation $\pi$ one associates a
collection of Hecke eigenvalues $\{ (z_1(\pi_x),\ldots,z_n(\pi_x))
\}_{x \in |X|-S}$, defined up to permutation. The strong multiplicity
one theorem due to Piatetski-Shapiro \cite{PS} says that this
collection determines $\pi$ up to an isomorphism.

\subsection{The Langlands correspondence and related results}

Now we can state the Langlands conjecture for $GL_n$ in the function
field case. It has been proved by Drinfeld \cite{Dr:icm,Dr:F} for
$n=2$ and by Lafforgue \cite{Laf} for $n>2$.

\begin{thm}    \label{langl}
There is a bijection between the sets ${\mc A}_n$ and ${\mc G}_n$
defined above which satisfies the following condition. If $\pi \in
{\mc A}_n$ corresponds to $\sigma \in {\mc G}_n$, then the sets of
points where they are unramified are the same, and for each $x$ from
this set we have
$$
(z_1(\pi_x),\ldots,z_n(\pi_x)) = (z_1(\sigma_x),\ldots,z_n(\sigma_x))
$$
up to permutation.
\end{thm}

In other words, if $\pi$ and $\sigma$ correspond to each other, then
the Hecke eigenvalues of $\pi$ coincide with the Frobenius eigenvalues
of $\sigma$ at all points where they are unramified.

\medskip

In addition to the Langlands correspondence, Drinfeld \cite{Dr:r} (for
$n=2$) and Lafforgue \cite{Laff,Laf} (for $n>2$) have also proved the
following result.

\begin{thm}    \label{rr}
At all points $x \in |X|$ where $\pi \in {\mc A}_n$ is
unramified, the Hecke eigenvalues $z_i(\pi_x)$ of $\pi$ are algebraic
numbers with (complex) absolute value equal to $1$.
\end{thm}

\begin{remark}
The statement saying that $|z_i(\pi_x)|=1$ is called the
Ramanujan-Petersson conjecture. Its classical analogue is the
Ramanujan conjecture, proved by Deligne, that the coefficients
$\tau(n)$ of the Ramanujan $\tau$--function
$$
q \prod_{m=1} (1-q^m)^{24} = \sum_{n=1}^\infty \tau(n) q^n
$$
satisfy the relation $|\tau(p)| \leq 2 p^{11/2}$ for a prime $p$. It
turns out that the numbers $\tau(p)$ may be written in the form
\begin{equation}    \label{alp}
\tau(p) = p^{11/2}(\al_p+\beta_p)
\end{equation}
where $\al_p$ and $\beta_p = \ol{\al}_p$ are the Hecke eigenvalues of
an automorphic representation of $GL_2(\AD_{\Q})$ corresponding to the
modular form $\tau(q)$ and therefore the Ramanujan conjecture is
equivalent to the statement that $|\al_p|=|\beta_p|=1$, which is
analogous to the condition appearing in \thmref{rr}.\qed
\end{remark}

\medskip

\thmref{rr} has a stunning corollary, known as the {\em Deligne purity
conjecture} \cite{De2}. Let $V$ be a normal algebraic variety over
$\Fq$. Then one defines, in the same way as above in the case of
curves, $\ell$--adic representations of the Galois group
$\on{Gal}(\ol{F}/F)$, where $F$ is the field of rational functions on
$V$ and $\ol{F}$ is its separable closure. The following theorem is
proved by reducing the statement to the case where $V$ is a curve and
applying Theorems \ref{langl} and \ref{rr}.

\begin{thm}
Let $\sigma$ be any irreducible $\ell$--adic representation
of $\on{Gal}(\ol{F}/F)$, which is everywhere unramified and has
determinant of finite order. Then the eigenvalues of the Frobenius
conjugacy classes at all closed points of $V$ are algebraic numbers
with (complex) absolute value $1$.
\end{thm}

\section{Elements of the proof of the Langlands conjecture}
\label{proof}

In this section we discuss some of the ideas and constructions
involved in the proof of the Langlands correspondence, \thmref{langl},
by Drinfeld and Lafforgue. In the first approximation, one can say
that the main idea is to realize this correspondence in the cohomology
of a certain moduli space of vector bundles on $X$ with some additional
structures, called ``shtukas'', and to use the
Grothendieck--Lefschetz formula to relate the traces of Hecke
correspondences acting on this cohomology with the numbers of fixed
points of these correspondences acting on the moduli space.

\subsection{From automorphic functions to vector bundles}    \label{vb}

Recall that automorphic representations of $GL_n(\AD)$ are realized in
the space ${\mc C}_{\on{cusp}}$ of functions on the quotient
$GL_n(F)\bs GL_n(\AD)$. Suppose that we are given an automorphic
representation $\pi$ of $GL_n(\AD)$ which is unramified at all points
of $X$. Then the space of $GL_n(\OO)$--invariants in $\pi$, where $\OO
= \prod_{x \in |X|} \OO_x$, is one-dimensional, spanned by the vector
$$v = \bigotimes_{x \in |X|} v_x,$$ where $v_x$ is defined in
\secref{aut repr}. Hence $v$ gives rise to a $GL_n(\OO)$--invariant
function on $GL_n(F)\bs GL_n(\AD)$, or equivalently, a function on the
double quotient $$GL_n(F)\bs GL_n(\AD)/GL_n(\OO).$$ The following key
observation is due to A. Weil.

\begin{lem}    \label{weil}
There is a bijection between the set $GL_n(F)\bs GL_n(\AD)/GL_n(\OO)$
and the set of isomorphism classes of rank $n$ vector bundles on $X$.
\end{lem}

\begin{proof}
Any rank $n$ bundle on $X$ may be trivialized on the formal disc $D_x
= \on{Spec} \OO_x$ around each point $x \in |X|$ and over $X - S$,
where $S$ is a sufficiently large subset of $|X|$, hence over the
generic point $X_{\on{gen}} = \on{Spec} F$ of $X$.

Let $B$ be the set of isomorphism classes of the data
$(\M,\vf_{\gen},(\vf_x))$, where $\M$ is a rank $n$ bundle on $X$, and
$\vf_{\gen}$ and $\vf_x, x \in |X|$, are the trivializations of $\M$
over the generic point of $X$ and the formal discs $D_x$,
respectively, i.e., isomorphisms $\vf_{\gen}: F^n \overset{\sim}\arr
\Gamma(X_{\on{gen}},\M)$, $\vf_x: \OO_x^n \overset{\sim}\arr
\Gamma(D_x,\M)$. The restrictions of $\vf_{\gen}$ and $\vf_x$ to the
punctured disc $D_x^\times = \on{Spec} \K_x$ then give us two
different trivializations of $\M|_{D^\times_x}$, which we denote by
the same symbols. Let $g_x = \vf_{\gen}^{-1} \circ \vf_x$ be the
corresponding transition function, which is an element of
$GL_n(\K_x)$.

Consider the map $b: B \arr GL_n(\AD)$ sending
$(\M,\vf_{\gen},(\vf_x))$ to $(g_x)_{x \in |X|} \in GL_n(\AD)$. It is
easy to see that this map is a bijection. Therefore there is bijection
between the set of isomorphism classes of rank $n$ bundles on $X$ and
the quotient of $GL_n(\AD)$ by the equivalence relations corresponding
to changes of the trivializations $\vf_{\on{gen}}$ and $\vf_x, x \in
|X|$. These equivalence relations amount to the left action of
$GL_n(F)$ and the right action of $GL_n(\OO)$ on $GL_n(\AD)$,
respectively. Hence we obtain the statement of the lemma.
\end{proof}

In order to apply the machinery of algebraic geometry one needs to
interpret sets like $GL_n(F)\bs GL_n(\AD)/GL_n(\OO)$ as sets of points
of algebraic varieties over $\Fq$. Once this is done, one can use
things like cohomology groups of varieties and techniques like the
Lefschetz fixed point formula. The above result gives precisely such
an interpretation of the double quotient $GL_n(F)\bs
GL_n(\AD)/GL_n(\OO)$. Namely, we obtain that it is the set of
$\Fq$--points of the moduli space $\on{Bun}_n$ of rank $n$ vector
bundles on $X$. To be precise, $\on{Bun}_n$ is not an algebraic
variety, but an algebraic {\em stack}, which means, roughly speaking,
that it looks locally like an algebraic variety quotiented out by the
action of an algebraic group (these actions are not free, and
therefore the quotient is no longer an algebraic variety), see
\cite{LMB} for the precise definition. But for our purposes this turns
out to be sufficient.

\subsection{Hecke correspondences}

Another important observation is the interpretation of the spherical
Hecke algebra in terms of {\em Hecke correspondences} in $\on{Bun}_n
\times \on{Bun}_n$. The spherical Hecke algebra ${\mc H} =
\bigotimes_{x \in |X|} {\mc H}_x$ naturally acts on the space of
functions on $GL_n(F)\bs GL_n(\AD)/GL_n(\OO)$, where it preserves the
one-dimensional subspaces $\pi^{GL_n(\OO)}$ of $GL_n(\OO)$--invariants
in the unramified irreducible representations $\pi$ of
$GL_n(\AD)$. Let us describe the algebra ${\mc H}_x$ in more
detail. According to formula \eqref{Hx}, ${\mc H}_x$ is isomorphic to
the algebra generated by the elementary symmetric functions in
$z_1,\ldots,z_n$, which we denote by $H_{1,x},\ldots,H_{n,x}$, and
$H_{n,x}^{-1}$. The action of the operator $H_{i,x}$ on the space of
functions on $GL_n(F)\bs GL_n(\AD)/GL_n(\OO)$ is given by the
following integral operators:
$$
\left( H_{i,x} \cdot f \right)(g) = \int_{M^i_n(\OO_x)} f(gh) dh,
$$
where $$M^i_n(\OO_x) = GL_n(\OO_x) \cdot \D_x^i \cdot GL_n(\OO_x)
\subset GL_n(F_x) \subset GL_n(\AD),$$ $\D_x^i$ is the diagonal matrix
whose first $i$ entries are equal to $t_x$ (a uniformizer at $x$), and
whose remaining $n-i$ entries are equal to $1$.

Now define the $i$th Hecke correspondence ${\mc H}ecke_i$ (in what
follows we will use the same notation for a vector bundle and for the
sheaf of its sections, which is a locally free coherent sheaf). It is
the moduli space of quadruples $$(\M,\M',x,\beta:
\M'\hookrightarrow\M),$$ where $\M',\M\in\Bun_n$, $x\in |X|$, and
$\beta$ is the embedding of the corresponding coherent sheaves of
sections $\beta:\M'\hookrightarrow\M$ such that $\M/\M'$ is supported
at $x$ and is isomorphic to the direct sum of $i$ copies of the
skyscraper sheaf $\OO_x = \OO_X/\OO_X(-x)$.

We have a correspondence
$$
\begin{array}{ccccc}
& & {\mc Hecke}_i & & \\
& \stackrel{\hl}\swarrow & & \stackrel{\supp\times\hr}\searrow & \\
\Bun_n & & & & X\times \Bun_n
\end{array}
$$
where $\hl(x,\M,\M')=\M$, $\hr(x,\M,\M')=\M'$, and $\supp(x,\M,\M') =
x$. We will use the same notation for the corresponding maps between
the sets of $\Fq$--points.

Let ${\mc H}ecke_{i,x} = \on{supp}^{-1}(x)$. This is a correspondence
in $\on{Bun}_n \times \on{Bun}_n$. Therefore it defines an operator on
the space of functions on $GL_n(F)\bs GL_n(\AD)/GL_n(\OO)$ which takes
a function $f$ to the function $\hr_!(\hl{}^*(f))$, where $\hr_!$ is
the operator of integration along the fibers of $\hr$. It is easy to
check that this operator is precisely the $i$th Hecke operator
$H_{i,x}$. Thus, we obtain an interpretation of the generators of the
spherical Hecke algebra ${\mc H}_x$ in terms of Hecke correspondences.

A general cuspidal automorphic representation $\pi$ is unramified away
from a finite set of points $S \subset |X|$, but for each $x \in S$
there exists a compact subgroup $K_x$ such that the space of
$K_x$--invariants in $\pi_x$ is non-zero. Without loss of generality
we may assume that $K_x$ is the congruence subgroup of $GL_n(\OO_x)$
whose elements are congruent to the identity modulo the $m_x$th power
of the maximal ideal of $\OO_x$. Consider the divisor $N = \sum_{x \in
S} m_x [x]$ on $X$. Denote by ${\mc A}_{n,N}$ the subset of ${\mc
A}_n$, which consists of the equivalence classes of those
representations $\pi$ which have a non-zero space of invariants with
respect to the compact subgroup $$K_N = \prod_{x \in S} K_x \times
\prod_{x \in |X|-S} GL_n(\OO_x).$$ If $\pi \in {\mc A}_{n,N}$, then
the space $\pi_N$ of $K_N$--invariants of $\pi$ embeds into the space
of functions on the double quotient $GL_n(F)\bs GL_n(\AD)/K_N$.

In a similar fashion to \lemref{weil} one identifies this double
quotient with the set of $\Fq$--points of the moduli stack
$\on{Bun}_{n,N}$ of rank $n$ vector bundles ${\mc M}$ on $X$ together
with the {\em level structure} at $N$, that is a trivialization of the
restriction of ${\mc M}$ to $N$, considered as a finite subscheme of
$X$. The Hecke algebra ${\mc H}^n_N$ of compactly supported functions
on $GL_n(\AD)$ bi-invariant with respect to $K_N$ acts on the space of
functions on $GL_n(F)\bs GL_n(\AD)/K_N$ preserving its subspace
$\pi^{K_N}$. This action may also be described in terms of Hecke
correspondences in $\on{Bun}_{n,N} \times \on{Bun}_{n,N}$ as in the
unramified case.

\subsection{Deligne's recurrence scheme}    \label{deligne}

The starting point of the proof of the Langlands correspondence is the
following recurrence scheme originally suggested by P. Deligne.

\begin{thm}    \label{recurrence}
Suppose that for each $n' = 1,\ldots,n-1$ we have a map $\rho_{n'}:
{\mc A}_{n'} \to {\mc G}_{n'}$ such that for each $\pi \in {\mc
A}_{n'}$ the Frobenius eigenvalues of $\rho_{n'}(\pi)$ coincide with
the Hecke eigenvalues of $\pi$ at all points where $\pi$ and
$\rho_{n'}(\pi)$ are unramified. Then there exists a map $\phi_n: {\mc
G}_n \to {\mc A}_n$ such that for each $\sigma \in {\mc G}_{n}$ the
Frobenius eigenvalues of $\sigma$ coincide with the Hecke eigenvalues
of $\phi_n(\sigma)$ at all points where $\sigma$ and
$\phi_{n}(\sigma)$ are unramified.
\end{thm}

The proof is based on several deep results: Grothendieck's functional
equation for the $L$--functions associated to $\ell$--adic Galois
representations \cite{Gr}, Laumon's product formula for the
$\ep$--constant appearing in this functional equation
\cite{Laumon:const}, and the ``converse theorems'' of Hecke, Weil and
Piatetski-Shapiro \cite{PS,CPS}. We refer the reader to
\cite{Laumon:const} and \cite{Laf} for more details.

The Langlands correspondence is known for $n=1$ by the abelian class
field theory. Therefore, in view of \thmref{recurrence}, in order to
establish the Langlands correspondence for all $n>1$ it is sufficient
to construct maps $\rho_{n}: {\mc A}_{n} \to {\mc G}_{n}$ satisfying
the conditions of \thmref{recurrence} for all $n>1$. How can this be
achieved?

A naive idea is to construct a natural representation of $GL_n(\AD)
\times G_F$ defined over $\Ql$ which decomposes (over $\oQl$) into a
direct sum
$$
\bigoplus_{\pi \in {\mc A}_n} \pi \otimes \sigma_\pi,
$$ where each $\sigma_\pi$ is an irreducible $n$--dimensional
$\ell$--adic representation of $G_F$. Then if the Hecke eigenvalues of
$\pi$ and Frobenius eigenvalues of $\sigma_\pi$ coincide, we can
construct the map $\rho_n$ by the formula $\pi \mapsto
\sigma_\pi$. Unfortunately, such a representation of $GL_n(\AD) \times
G_F$ does not exist (this is explained in \cite{K}). Instead, Drinfeld
proposed to construct a representation of the product $GL_n(\AD)
\times G_F \times G_F$ which decomposes as
$$
\bigoplus_{\pi \in {\mc A}_n} \pi \otimes \sigma_\pi \otimes
\sigma_\pi^\vee,
$$
where $\sigma_\pi^\vee$ denotes the representation contragredient to
$\sigma_\pi$ (here and below we will be ignoring the Tate twists).

Drinfeld (for $n=2$) and Lafforgue (for $n>2$) have constructed such a
representation in the ``essential'' part of the direct limit of the
middle $\ell$--adic cohomologies of the moduli spaces of {\em shtukas}
with level structures. In the next section we will introduce these
moduli spaces.

\subsection{Moduli spaces of shtukas}

As is well-known by now, ``shtuka'' is a Russian word that may be
loosely translated as a ``widget''. This term was used by Drinfeld for
the following objects that he introduced.

First note that the Frobenius endomorphism $\on{Fr}$ of $\Fq$ defined
by the formula $y \mapsto y^q$ induces a map $\on{Id} \times \on{Fr}:
\ol{X} \to \ol{X}$, where $\ol{X} = X \underset{\on{Spec} \Fq}\times
\on{Spec} \oFq$. Given a vector bundle $\E$ on $\ol{X}$, we denote by
$^\tau\E$ the vector bundle $(\on{Id} \times \on{Fr})^*(\E)$.

Recall that the skyscraper sheaf supported at $z \in \ol{X}$ is the
coherent sheaf whose stalk at $y \in \ol{X}$ is one-dimensional if
$y=z$ and is equal to $0$ if $y \neq z$.

\begin{definition}
{\em A shtuka of rank $n$ on $\ol{X}$ is a vector bundle $\E$ of rank $n$
on $\ol{X}$ together with a diagram
$$
\E \overset{j}\hookrightarrow \E' \overset{t}\hookleftarrow {}^\tau\E
$$
where $\E'$ is another rank $n$ vector bundle on $\ol{X}$ and $t,j$
are injections of the corresponding sheaves of sections such that
their cokernels are the skyscraper sheaves supported at the points
$0,\infty \in \ol{X}$, called the zero and the pole of the shtuka,
respectively.}
\end{definition}

Note that a vector bundle $\E$ on $\ol{X}$ equipped with an
isomorphism $^\tau\E \simeq \E$ is the same thing as a vector bundle
on $X$. So a shtuka is a mild generalization of the notion of vector
bundle on $X$. Indeed, for any shtuka the bundles $^\tau\E$ and $\E$
are isomorphic over $\ol{X} - \{ 0,\infty \}$ and they differ in the
simplest possible way at the points $0$ and $\infty$.

The reader may find more information on shtukas and closely related
objects, Drinfeld modules, including explicit examples, in \cite{Go}.

Let $N = \on{Spec} \OO_N$ be a finite subscheme of $X$ (equivalently,
a divisor on $X$). A shtuka {\em with level structure} $N$ is a shtuka
$\wt\E = (\E,\E',j,t)$ such that the points $0$ and $\infty$ avoid
$N$, together with a trivialization of the restriction of $\E$ to $N$,
i.e., an isomorphism $\E \underset{\OO_{\ol{X}}}\otimes \OO_{\ol{N}}
\simeq \OO_{\ol{N}}^{\oplus n}$ so that the induced trivialization of
the restriction of $^\tau\E$ to $N$ is compatible with the isomorphism
$$
t \circ j: \E|_{\ol{X}\bs\{0,\infty\}} \overset{\sim}\to
{}^\tau\E|_{\ol{X}\bs\{0,\infty\}}.
$$

One defines similarly shtukas over $X \times S$, where $S$ is any
scheme over $\Fq$. This enables one to define an algebraic stack
$\on{Cht}^n_N$ classifying shtukas of rank $n$ with level structure at
$N$ on $X$. Drinfeld proved that this is a smooth Deligne-Mumford
stack, which means, roughly, that locally it looks like the quotient
of a smooth algebraic variety by the action of a finite group. It is
equipped with a natural morphism to $(X-N) \times (X-N)$ (taking the
pole and zero of a shtuka) of relative dimension $2n-2$. It carries
the action of the ``partial'' Frobenius endomorphisms $\on{Frob}_0$
and $\on{Frob}_\infty$ corresponding to the zero and pole of the
shtuka. In the same way as in \secref{vb} one defines Hecke
correspondences in $\on{Cht}^n_N \underset{(|X|-N)^2}\times
\on{Cht}^n_N$ which realize the Hecke algebra ${\mc H}^n_N$ of
compactly supported functions on $GL_n(\AD)$ bi-invariant with respect
to $K_N$.

In addition, one has a natural action on $\on{Cht}^n_N$ of the Picard
group of line bundles on $X$ (by tensoring with $\E$ and $\E'$).
\lemref{weil} implies that the Picard group is isomorphic to the
quotient $F^\times\bs \AD^\times/\OO^\times$. Let us pick an element
$a$ of degree $1$ in $\AD^\times$. Lafforgue denotes the quotient of
$\on{Cht}^n_N$ by the action of the cyclic group generated by the
corresponding line bundle by $\on{Cht}^n_N/a^{\Z}$. This algebraic
stack still carries the above actions of $\on{Frob}_0,
\on{Frob}_\infty$ and ${\mc H}^n_N$.

\subsection{Strategy of the proof}

Let $q', q''$ be the two projections $X^2 \to X$. For each
representation $\sigma$ of $G_F$ we obtain by pull-back two
representations, $q'{}^*(\sigma)$ and $q''{}^*(\sigma)$ of the Galois
group $G_{\wt{F}}$ of the field $\wt{F}$ of functions on
$X^2$. Consider the $\ell$--adic cohomology with compact support of
$\on{Cht}^n_N/a^{\Z}$, over the generic point of $X \times X$. This is
naturally a representation of the Hecke algebra ${\mc H}^n_N$ and of
$G_{\wt{F}}$, whose actions commute with each other. We would like to
isolate in this cohomology a subspace that decomposes as
\begin{equation}    \label{decomp}
\bigoplus_{\pi \in {\mc A}^a_{n,N}} \pi_{N} \otimes q'{}^*(\sigma_\pi)
\otimes q''{}^*(\sigma_\pi^\vee),
\end{equation}
where $\pi_N$ is the space of $K_N$--invariants in $\pi$, each
$\sigma_\pi$ is an irreducible $n$--dimensional $\ell$--adic
representation of $G_F$ unramified in $|X|-N$, and ${\mc A}^a_{n,N}$
is a subset of ${\mc A}_{n,N}$, which consists of those
representations on which $a \in GL_n(\AD)$ acts as the identity. If we
could show that the Hecke eigenvalues of each $\pi$ in the above
formula coincide with the Frobenius eigenvalues of $\sigma_\pi$, then
in view of the discussion at the end of \secref{deligne}, this would
prove the Langlands correspondence.

Suppose we could isolate this subspace in the cohomology. Then we
would need to compute the traces of the operators of the form $f
\times (\on{Frob}_0)^s \times (\on{Frob}_\infty)^s$ on this subspace,
where $f \in {\mc H}^n_N$ and $0,\infty \in |X-N|$. If we could
establish that these traces are equal to
\begin{multline}    \label{lef}
\sum_{\pi \in {\mc A}^a_{n,N}} \on{Tr}_{\pi_N}(f)
\left(z_1(\pi_\infty)^{-s/\deg(\infty)} + \ldots +
z_n(\pi_\infty)^{-s/\deg(\infty)} \right)\times \\
\left(z_1(\pi_0)^{s/\deg(0)} + \ldots + z_n(\pi_0)^{s/\deg(0)}
\right),
\end{multline}
then it would not be difficult to prove that our space is indeed
isomorphic to \eqref{decomp} with all the required
compatibilities. But how could we possibly identify the trace over the
cohomology with the sum \eqref{lef}?

First of all, we need to apply the {\em Grothendieck-Lefschetz} fixed
point formula. It expresses the alternating sum of the traces of
correspondences over $\ell$--adic cohomologies of a smooth variety as
the number of the fixed points of these correspondences. Suppose we
could apply this formula to the correspondences $f \times
(\on{Frob}_0)^s \times (\on{Frob}_\infty)^s$ acting on
$\on{Cht}^n_N/a^{\Z}$. We should then compare the number of fixed
points appearing in this formula with the {\em Arthur-Selberg trace
formula}. This formula describes the traces of the operators like $f
\in {\mc H}^n_N$ on the space ${\mc C}_{\on{cusp}} = \bigoplus_{\pi
\in {\mc A}_n} \pi$ of cuspidal automorphic functions in terms of
orbital integrals in $GL_n(\AD)$. We may hope to relate these orbital
integrals to the numbers of fixed points in the moduli spaces of
shtukas. That would give us the desired expression for the trace of
our correspondences acting on the $\ell$--adic cohomology as the sum
\eqref{lef}.

\subsection{From a dream to reality}

This is the general strategy of Drinfeld and Lafforgue. Unfortunately,
literally it cannot work, because $\on{Cht}^n_N/a^{\Z}$ is not
(quasi)compact: it cannot be covered by finitely many open subsets of
finite type. For this reason its cohomology is infinite-dimensional
and there are infintely many fixed points, and so one cannot apply the
Grothendieck-Lefschetz formula. To remedy this, Lafforgue introduces
open substacks $\on{Cht}^{n,p}_N/a^{\Z}$ of finite type in
$\on{Cht}^{n,p}_N/a^{\Z}$. They are labeled by ``Harder-Narasimhan
polygons'' $p: [0,n] \to {\mathbb R}_+$. Furthermore, using the
Arthur-Selberg trace formula (see \cite{Laff}) Lafforgue computes the
numbers of fixed points of the correspondences $f \times
(\on{Frob}_0)^s \times (\on{Frob}_\infty)^s$ in the set of
$\Fq$--points of $\on{Cht}^{n,p}_N/a^{\Z}$. The answer is the
expression \eqref{lef} plus the the sum of terms that correspond to
representations of $GL_n(\AD)$ induced from parabolic subgroups (i.e.,
those which are not cuspidal).

The problem however is that (except for the case when it is equal to
the identity) the correspondence $f$ does not stabilize the open
subset $\on{Cht}^{n,p}_N/a^{\Z}$. Therefore we cannot interpret the
number of fixed points as the trace over the cohomology of
$\on{Cht}^{n,p}_N/a^{\Z}$. To fix this problem, Lafforgue (and
Drinfeld for $n=2$) introduced compactifications
$\ol{\on{Cht}^{n,p}_N}/a^{\Z}$ of $\on{Cht}^{n,p}_N/a^{\Z}$ by
allowing certain degenerations of shtukas. The Hecke correspondences
can now be extended to $\ol{\on{Cht}^{n,p}_N}/a^{\Z}$. Unfortunately,
the compactifications $\ol{\on{Cht}^{n,p}_N}/a^{\Z}$ are singular
(unless $N=\emptyset$). Therefore a priori these Hecke correspondences
do not induce linear operators on the cohomologies of
$\ol{\on{Cht}^{n,p}_N}/a^{\Z}$.

In the case when $n=2$, Drinfeld has shown that
$\ol{\on{Cht}^{2,p}_N}/a^{\Z}$ is quasi-smooth, i.e., its cohomology
exhibits Poincar\'e duality. Therefore the Hecke correspondences act
on the cohomology of $\ol{\on{Cht}^{2,p}_N}/a^{\Z}$, and Drinfeld was
able to relate their traces to the numbers of fixed points. This
allowed him to prove that the cuspidal part of the middle (second in
this case) cohomology of $\ol{\on{Cht}^{2,p}_N}/a^{\Z}$ indeed
decomposes according to formula \eqref{decomp} and hence realizes the
Langlands correspondence for $GL_2$.

For $n>2$ the stacks $\ol{\on{Cht}^{n,p}_N}/a^{\Z}$ are no longer
quasi-smooth in general. To deal with this problem Lafforgue
introduced another open subset $\ol{\on{Cht}^{n,p}_N}'/a^{\Z}$ of
$\ol{\on{Cht}^{n,p}_N}/a^{\Z}$ defined by the condition that the
degenerations of the shtuka avoid $N$. The stack
$\ol{\on{Cht}^{n,p}_N}'/a^{\Z}$ is smooth and its complement in
$\ol{\on{Cht}^{n,p}_N}/a^{\Z}$ is a divisor with normal
crossings. Moreover, it turns out that the (normalized) Hecke
correspondences defined in $\ol{\on{Cht}^{n,p}_N}'/a^{\Z}$ stabilize
$\ol{\on{Cht}^{n,p}_N}'/a^{\Z}$.

Thus, Lafforgue ended up with three different objects, each with some
``good'' and ``bad'' properties. Indeed, $\on{Cht}^n_N/a^{\Z}$ carries
an action of the Hecke algebra by correspondences, but is not of
finite type. The stacks $\on{Cht}^{n,p}_N/a^{\Z}$ are of finite type
and it is possible to compute the numbers of fixed points of the Hecke
correspondences there, but the correspondences themselves do not
preserve them. Finally, on $\ol{\on{Cht}^{n,p}_N}'/a^{\Z}$ one can
write down the Grothendieck-Lefschetz formula for the Hecke
correspondences, but it is not clear how to compute either side of
this formula.

Lafforgue's ingenious trick is to separate inside the cohomology with
compact support of all three objects the ``essential'' part and the
``negligible'' part and to show that the ``essential'' part is the
same for all of them. Namely, the negligible part consists of those
representations of $G_{\wt{F}}$ which appear as direct factors in the
tensor products of the form $q'{}^*(\sigma') \otimes
q''{}^*(\sigma'')$, where $\sigma'$ and $\sigma''$ are
$G_F$--representations of dimension less than $n$. The rest is the
essential part.

In order to compute this essential part, Lafforgue needed to overcome
some formidable technical difficulties. The original
Grothendieck-Lefschetz formula is applicable when we have a smooth
proper algebraic variety. Here $\ol{\on{Cht}^{n,p}_N}'/a^{\Z}$ is not
a variety but an algebraic stack which is not even proper. Building
upon earlier works of Pink and Fujiwara on the Deligne conjecture,
Lafforgue proved a new version of the Grothendieck-Lefschetz formula
which enabled him to express the trace of a Hecke correspondence and a
power of the Frobenius endomorphism on the cohomology of
$\ol{\on{Cht}^{n,p}_N}'/a^{\Z}$ as the number of fixed points in
$\on{Cht}^{n,p}_N/a^{\Z}$ plus a sum of terms corresponding to various
boundary strata in the complement of $\on{Cht}^{n,p}_N/a^{\Z}$. By a
complicated recurrence argument on $n$ he showed that the latter are
all negligible. Therefore he was able to identify the traces of the
Hecke operators and the Frobenius endomorphisms acting on the
essential part of the middle cohomology with compact support of
$\on{Cht}^{n,p}_N/a^{\Z}$ (which coincides with that of
$\ol{\on{Cht}^{n,p}_N}'/a^{\Z}$) with formula \eqref{lef}. This
completed his proof of the Langlands correspondence for $GL_n$.

In addition, he proved \thmref{rr}, because the Frobenius eigenvalues
of irreducible representations of $G_F$ are now realized as traces of
Frobenius endomorphisms acting on $\ell$--adic cohomology, so that one
can apply Deligne's results on Weil's conjectures \cite{De1,De2}.

\section{The geometric Langlands conjecture}    \label{geometric}

The geometric reformulation of the Langlands conjecture allows one to
state it for curves defined over an arbitrary field, not just a finite
field. For instance, it may be stated for complex curves, and in this
setting one can apply methods of complex algebraic geometry which are
unavailable over finite fields. Hopefully, this will eventually help
us understand better the general underlying patterns of the Langlands
correspondence. In this section we will formulate the geometric
Langlands conjecture for $GL_n$ and discuss briefly its recent proof
due to D. Gaitsgory, K. Vilonen and the author.

\subsection{Galois representations as local systems}

What needs to be done to reformulate the Langlands conjecture
geometrically?

As we indicated in \secref{galois}, the Galois group $G_F$ should be
viewed as a kind of fundamental group, and so its $\ell$--adic
representations unramified away from a divisor $N$ should be viewed as
local systems on $X-N$. This is indeed possible if one defines local
systems in terms of $\ell$--adic sheaves.

Let us discuss them briefly. Let $V$ be an algebraic variety over
$\Fq$, and $\ell$ a prime which does not divide $q$. Then one defines
the category of $\ell$--adic sheaves on $V$. The construction involves
several steps (see, e.g., \cite{Milne,Weil}). First one considers
locally constant $\Z/\ell^m \Z$--sheaves on $V$ in the \'etale
topology (in which the role of open subsets is played by \'etale
morphisms $U \to V$). A $\Z_\ell$--sheaf on $V$ is by definition a
system $(\F_m)$ of locally constant $\Z/\ell^m \Z$--sheaves satisfying
natural compatibilities. Then one defines the category of
$\Ql$--sheaves by killing the torsion sheaves in the category of
$\Z_\ell$--sheaves. In a similar fashion one defines the category of
$E$--sheaves on $V$, where $E$ is a finite extension of
$\Ql$. Finally, one takes the direct limit of the categories of
$E$--sheaves on $X$, and the objects of this category are called the
lisse $\ell$--adic sheaves on $V$, or $\ell$--{\em adic local
systems}. An $\ell$--adic local system on $X-N$ of rank $n$ is the
same as an $n$--dimensional $\ell$--adic representation of $G_F$
unramified everywhere on $X-N$.

Thus, we can now interpret Galois representations geometrically: these
are $\ell$--adic local systems on $X$. The notion of local system
makes sense if $X$ is defined over other fields. For example, if $X$
is a smooth projective curve over $\C$, a local system is a vector
bundle on $X$ with a flat connection, or equivalently, a homomorphism
from the fundamental group $\pi_1(X)$ to $GL_n(\C)$.

\subsection{From functions to sheaves}    \label{grot}

Next, we wish to interpret geometrically automorphic
representations. Let us restrict ourselves to unramified
representations. As we explained in \secref{vb}, one can attach to
such a representation $\pi$ a non-zero function $f_\pi$ (unique up to
a scalar) on $GL_n(F)\bs GL_n(\AD)/GL_n(\OO)$, which is an
eigenfunction of the Hecke algebras ${\mc H}_x, x \in |X|$. In fact,
this function completely determines the representation $\pi$, so
instead of considering the set of equivalence classes of unramified
cuspidal representations of $GL_n(\AD)$, one may consider the set of
unramified automorphic functions associated to them (defined up to a
scalar).

The key step in the geometric reformulation of this notion is the
Grothendieck {\em fonctions--faisceaux} dictionary. Let $V$ be an
algebraic variety over $\Fq$. One generalizes the above definition of
an $\ell$--adic local system on $V$ by allowing the
$\Z/\ell^n\Z$--sheaves ${\mc F}_n$ to be constructible, i.e., for
which there exists a stratification of $V$ by locally closed
subvarieties $V_i$ such that the sheaves ${\mc F}|_{V_i}$ are locally
constant. As a result, one obtains the notion of a constructible
$\ell$--adic sheaf on $V$, or an $\ell$--adic sheaf, for brevity. Let
${\mc F}$ be such a sheaf and $x$ be an ${\mathbb F}_{q_1}$--point of
$V$, where $q_1=q^m$. Then one has the Frobenius conjugacy class
$\on{Fr}_x$ acting on the stalk ${\mc F}_x$ of ${\mc F}$ at $x$. Hence
we can define a function $\text{\tt f}_{q_1}({\mc F})$ on the set of
${\mathbb F}_{q_1}$--points of $V$, whose value at $x$ is
$\on{Tr}(\on{Fr}_x,{\mc F}_x)$.

More generally, if $\K$ is a complex of $\ell$--adic sheaves, one
defines a function $\text{\tt f}_{q_1}({\mc K})$ on $V({\mathbb
F}_{q_1})$ by taking the alternating sums of the traces of $\on{Fr}_x$
on the stalk cohomologies of $\K$ at $x$. The map $\K \to \text{\tt
f}_{q_1}({\mc K})$ intertwines the natural operations on sheaves with
natural operations on functions (see \cite{Laumon:const},
Sect. 1.2). For example, pull-back of a sheaf corresponds to the
pull-back of a function, and push-forward of a sheaf {\em with compact
support} corresponds to the fiberwise integration of a function (this
follows from the Grothendieck-Lefschetz trace formula).

We wish to identify a natural abelian category in the derived category
of $\ell$--adic sheaves such that the map ${\mc K} \mapsto (\text{\tt
f}_{q_1}({\mc K}))_{q_1=q^m}$ is injective. The naive category of
$\ell$--adic sheaves is not a good choice for various reasons; for
instance, it is not stable under the Verdier duality. The correct
choice is the abelian category of {\em perverse sheaves}. These are
complexes of $\ell$--adic sheaves on $V$ satisfying certain
restrictions on the degrees of their non-zero stalk cohomologies (see
\cite{BBD}). Examples are $\ell$--adic local systems on a smooth
variety $V$, placed in the cohomological degree equal to $-\dim
V$. Unlike ordinary sheaves, the perverse sheaves have the following
remarkable property: an irreducible perverse sheaf on a variety $V$ is
completely determined by its restriction to an arbitrary open dense
subset, provided that this restriction is non-zero.

Optimistically, one can hope that all ``interesting'' functions on
$V({\mathbb F}_{q_1})$ come from perverse sheaves by taking the traces
of the Frobeniuses (or is it Frobenii?). Unramified automorphic
functions on $GL_n(F)\bs GL_n(\AD)/GL_n(\OO)$ are certainly
``interesting''. Therefore one hopes that they all come from perverse
sheaves on the moduli stack $\on{Bun}_n$, whose set of points is
$GL_n(F)\bs GL_n(\AD)/GL_n(\OO)$, according to \lemref{weil}.

Thus, we have identified the geometric objects which should replace
unramified automorphic functions: these are perverse sheaves on
$\on{Bun}_n$. This concept also makes sense if $X$ is defined over
other fields, for example, the field of complex numbers (see, e.g.,
\cite{GM}).

\subsection{Hecke eigensheaves and the geometric Langlands conjecture}
\label{Hecke functors}

But how to formulate the Hecke eigenfunction condition which
unramified automorphic functions satisfy in sheaf-theoretic terms? The
key is the description of the Hecke operators in terms of the Hecke
correspondences that was explained in \secref{vb}. We use these
correspondences to define the {\em Hecke functors} $\on{H}_i$ from the
category of perverse sheaves on $\on{Bun}_n$ to the derived category
of sheaves on $X \times \on{Bun}_n$ by the formula (here and below we
will ignore cohomological shifts and Tate twists)
\begin{equation}    \label{formula H1}
\He_i(\K) = (\supp\times\hr)_! \hl{}^*(\K).
\end{equation}

\begin{definition}
{\em Consider an $\ell$--adic local system $E$ of rank $n$ on $X$. A
perverse sheaf $\K$ on $\Bun_n$ is called a {\em Hecke eigensheaf}
with respect to $E$, if $\K\neq 0$ and we are given isomorphisms}
\begin{equation} \label{eigen-property}
\He^i_n(\K)\simeq \wedge^i E \boxtimes \K, \quad i=1,\ldots,n,
\end{equation}
{\em where $\wedge^i E$ is the $i$th exterior power of $E$}.
\end{definition}

Let $\sigma$ be an $n$--dimensional unramified $\ell$--adic
representation of $G_F$ and $E_\sigma$ the corresponding $\ell$--adic
local system on $X$. Then $$\on{Tr}(\on{Fr}_x,E_x) =
\on{Tr}(\sigma(\on{Fr}_x),\oQl^n) = \sum_{i=1}^n z_i(\sigma_x)$$ (see
\secref{galois} for the definition of $z_i(\sigma_x)$), and so
$$
\on{Tr}(\on{Fr}_x,\wedge^i E_x) =
s_i(z_1(\sigma_x),\ldots,z_n(\sigma_x)),
$$
where $s_i$ is the $i$th elementary symmetric polynomial.
Therefore we find that the function
$\text{\tt f}_{q}(\K)$ on $GL_n(F)\bs GL_n(\AD)/GL_n(\OO)$ associated
to a Hecke eigensheaf $\K$ satisfies
$$
H_{i,x} \cdot \text{\tt f}_{q}(\K) =
s_i(z_1(\sigma_x),\ldots,z_n(\sigma_x)) \text{\tt f}_{q}(\K)
$$
(up to some $q$--factors). In other words, $\text{\tt f}_{q}(\K)$ is a
Hecke eigenfunction whose Hecke eigenvalues are equal to the Frobenius
eigenvalues of $\sigma$. Hence we are naturally led to the following
geometric Langlands conjecture, which is due to Drinfeld and Laumon
\cite{Laumon:cor1}.

The statement of this conjecture was proved by Deligne for $GL_1$ (we
recall it in the next section) and by Drinfeld in the case of $GL_2$
\cite{Dr1}. These works motivated the conjecture in the case of
$GL_n$, which has recently been proved in \cite{FGV,Ga}. So we now
have the following result, in which $X$ is a smooth projective
connected curve defined either over a finite field or over the field
of complex numbers. As before, we denote by $\Bun_n$ the moduli stack
of rank $n$ bundles on $X$. It is a disjoint union of connected
components $\Bun^d_n$ corresponding to vector bundles of degree $d$.

\begin{thm}    \label{glc}
For each irreducible rank $n$ local system $E$ on $X$ there exists a
perverse sheaf $\Aut_E$ on $\Bun_n$, irreducible on each connected
component $\Bun^d_n$, which is a Hecke eigensheaf with respect to $E$.
\end{thm}

\subsection{Geometric abelian class field theory}    \label{gacft}

Let us consider Deligne's proof of the $n=1$ case of \thmref{glc} (see
\cite{Laumon:cor1}). In this case $Bun^d_n$ is essentially the
component $\on{Pic}_d$ of the Picard variety $\on{Pic}$ of $X$
classifying the line bundles on $X$ of degree $d$. In order to prove
\conjref{glc} in this case, we need to assign to each rank one local
system $E$ on $X$ a perverse sheaf $\Aut_E$ on $\on{Pic}$ which
satisfies the following Hecke eigensheaf property:
\begin{equation}    \label{hecke one}
\hl{}^*(\Aut_E) \simeq E \boxtimes \Aut_E,
\end{equation}
where $\hl: X \times \on{Pic} \to \on{Pic}$ is given by $(\Ll,x) \mapsto
\Ll(x)$.

Consider the Abel-Jacobi map $\pi_d: S^d X \to \on{Pic}_d$ sending the
divisor $D$ to the line bundle $\OO(D)$. If $d>2g-2$, then $\pi_d$ is
a projective bundle, with the fibers $\pi_d^{-1}(\Ll) = {\mathbb P}
H^0(X,\Ll)$. It is easy to construct a local system $E^{(d)}$ on
$\bigcup_{d>2g-2} S^d X$ satisfying an analogue of the Hecke
eigensheaf property
\begin{equation}    \label{hecke rank one}
\wt{h}^\leftarrow{}^*(E^{(d)}) \simeq E \boxtimes E^{(d)},
\end{equation}
where $\wt{h}^\leftarrow: S^d X \times X \to S^{d+1} X$ is given by
$(D,x) \mapsto D+[x]$. Namely, let $$\on{sym}^d: X^n \to S^n X$$ be
the symmetrization map and set $E^{(d)} = (\on{sym}^d_*(E^{\boxtimes
n}))^{S_d}$.

So we have local systems $E^{(d)}$ on $S^d X, d>2g-2$, which form a
Hecke eigensheaf, and we need to prove that they descend to
$\on{Pic}_d$ under the Abel-Jacobi map $\pi_d$. In other words, we
need to prove that the restriction of $E^{(d)}$ to each fiber of
$\pi_d$ is a constant sheaf. Since $E^{(d)}$ is a local system, these
restrictions are locally constant. But the fibers of $\pi_d$ are
projective spaces, hence simply-connected. Therefore any locally
constant sheaf along the fiber is constant! So there exists a local
system $\Aut^d_E$ on $\on{Pic}_d$ such that $E^{(d)} =
\pi_d^*(\Aut^d_E)$. Formula \eqref{hecke rank one} implies that the
sheaves $\Aut^d_E$ form a Hecke eigensheaf on $\bigcup_{d>2g-2}
\on{Pic}_d$. We extend them to the remaining components of $\on{Pic}$
by using the Hecke eigensheaf property \eqref{hecke one}. Namely, we
pick a point $x \in |X|$ and set $\Aut_E^d = E_x^* \otimes
\hl_x{}^*(\Aut_E^{d+1})$, where $\hl_x(\Ll) = \Ll(x)$. Then the fact
that the restrictions of $E^{(d)}$ to the fibers of $\pi_d$ is
constant implies that the resulting sheaves $\Aut^d_E$ for $d\leq
2g-2$ do not depend on the choice of $x$. Thus, we obtain a Hecke
eigensheaf on the entire $\on{Pic}$, and this completes Deligne's
proof of the geometric Langlands conjecture for $n=1$.

Let us consider the case when $X$ is defined over a finite field. Then
to the sheaf $\Aut_E$ we attach a function on $F^\times\bs
\AD^\times/\OO^\times$, which is the set of $\Fq$--points of
$\on{Pic}$. This function is a Hecke eigenfunction $f_\sigma$ with
respect to a one-dimensional Galois representation $\sigma$
corresponding to $E$, i.e., it satisfies the equation
$f_\sigma(\Ll(x)) = \sigma(\on{Fr}_x) f_\sigma(\Ll)$ (since $\sigma$
is one-dimensional, we do not need to take the trace). We could try to
construct this function proceeding in the same way as above. Namely,
we define first a function $f'_\sigma$ on the set of all divisors
on $X$ by the formula
$$
f'_\sigma\left(\sum_i n_i [x_i] \right) = \prod_i
\sigma(\on{Fr}_{x_i})^{n_i}.
$$
This function clearly satisfies an analogue of the Hecke eigenfunction
condition. It remains to show that the function $f'_\sigma$ descends
to $\on{Pic}(\Fq)$, namely, that if two divisors $D$ and $D'$ are
rationally equivalent, then $f'_\sigma(D) = f'_\sigma(D')$. This is
equivalent to the identity
$$
\prod_i \sigma(\on{Fr}_{x_i})^{n_i} = 1, \qquad \on{if} \quad \sum_i
n_i [x_i] = (g),
$$
where $g$ is an arbitrary rational function on $X$. This identity is a
non-trivial reciprocity law which has been proved in the abelian class
field theory, by Lang and Rosenlicht (see \cite{Serre}).

It is instructive to contrast this to Deligne's geometric proof
reproduced above. When we replace functions by sheaves we can use
additional information which is ``invisible'' at the level of
functions, such as the fact that that the sheaf corresponding to the
function $f'_\sigma$ is locally constant and that the fibers of the
Abel-Jacobi map are simply-connected. This is one of the main
motivations for studying the Langlands correspondence in the geometric
setting.

\subsection{The idea of the proof for general $n$}

Observe that for large $d$ the variety $S^d X$ may be interpreted as
the moduli space $\on{Bun}'_1$ of pairs $(\Ll,s)$, where $\Ll$ is a
line bundle on $X$ and $s$ is its section. We first constructed a
Hecke eigensheaf on $\on{Bun}'_1$ and then showed that it descends to
$\on{Bun}_1$. This is the main idea of the construction of $\Aut_E$
for general $n$ as well.

At the level of functions this construction is due to Weil \cite{W}
and Jacquet-Langlands \cite{JL} for $n=2$, and Shalika \cite{Sh} and
Piatetski-Shapiro \cite{PS} for general $n$. They attach to an
unramified $n$--dimensional representation $\sigma$ of $G_F$, a
function $f'_\sigma$ on the set of isomorphism classes of pairs
$(\M,s)$, where $\M \in \Bun_n$ is a rank $n$ bundle on $X$ and $s$ is
a regular non-zero section of $\M$. Then it remains to show that this
function is independent of the section, i.e., descends to the the set
of isomorphism classes of rank $n$ bundles on $X$, which is the double
quotient $GL_n(F)\bs GL_n(\AD)/GL_n(\OO)$.

We reformulate this geometrically. Let $\Bun'_n$ be the moduli stack
of pairs $(\M,s)$, where $\M$ is a rank $n$ bundle on $X$ and $s$ is a
regular non-zero section of $L$. Let $E$ be an irreducible rank $n$
local system on $X$. Building on the ideas of Drinfeld's work
\cite{Dr1}, Laumon gave a conjectural construction of the Hecke
eigensheaf $\Aut_E$ in \cite{Laumon:cor1,Laumon:cor2}. More precisely,
he attached to each rank $n$ local system $E$ on $X$ a complex of
perverse sheaves $\Aut'_E$ on $\Bun'_n$ and conjectured that if $E$ is
irreducible then this sheaf descends to a perverse sheaf $\Aut_E$ on
$\Bun_n$ (irreducible on each component), which is a Hecke eigensheaf
with respect to $E$.

In the paper \cite{FGKV} it was shown that the function on the set of
points $\Bun'_n(\Fq)$ associated to $\Aut'_E$ agrees with the function
$f'_\sigma$ constructed previously, as anticipated by Laumon
\cite{Laumon:cor2}. This provided a consistency check for Laumon's
construction. Next, in \cite{FGV}, Gaitsgory, Vilonen and myself
formulated a certain {\em vanishing conjecture} and proved that
Laumon's construction indeed produces a perverse sheaf $\Aut_E$ on
$\Bun_n$ with the desired properties whenever the vanishing conjecture
holds for $E$.  In other words, the vanishing conjecture implies the
geometric Langlands conjecture, for curves over any field. Moreover,
in the case when this field is $\Fq$, we derived the vanishing
conjecture (and hence the geometric Langlands conjecture) from the
results of Lafforgue \cite{Laf}. Finally, Gaitsgory \cite{Ga} gave
another proof of the vanishing conjecture, valid for curves both over
$\Fq$ and over $\C$. In the next section we will state this vanishing
conjecture.

\subsection{The vanishing conjecture}

Denote by $\Coh_n$ the stack classifying coherent shea\-ves on $X$ of
generic rank $n$, and by $\Coh_n^d$ its connected component
corresponding to coherent shea\-ves of degree $d$.

In \cite{Laumon:cor1} Laumon associated to an arbitrary local system
$E$ of rank $n$ on $X$ a perverse sheaf ${\mathcal L}_E$ on
$\Coh_0$. Let us recall his construction. Denote by
$\Coh_0^{\on{rss}}$ the open substack of $\Coh_0$ corresponding to
regular semisimple torsion sheaves. Thus, a geometric point of
$\Coh_0$ belongs to $\Coh_0^{\on{rss}}$ if the corresponding coherent
sheaf on $X$ is a direct sum of skyscraper sheaves of length one
supported at distinct points of $X$. Let
$\Coh_0^{\on{rss},d}=\Coh_0^{\on{rss}}\cap \Coh_0^d$. We have a
natural smooth map $(S^d X-\Delta)\to \Coh_0^{\on{rss},d}$.

Let $E^{(d)}$ be the $d$th symmetric power of $E$ defined as in
\secref{gacft}, i.e., $E^{(d)} = (\on{sym}_*(E^{\boxtimes d}))^{S_d}$,
where $\on{sym}: X^d \to S^d X$. This is a perverse sheaf on $S^d X$,
though it is not a local system if $n>1$. However, its restriction
$E^{(d)}|_{S^d X-\Delta}$ is a local system, which descends to a local
system $\ovc\Ll{}_E^d$ on $\Coh_0^{\on{rss},d}$. The perverse sheaf
$\Ll^d_E$ on $\Coh^d_0$ is by definition the canonical extension of
$\ovc\Ll{}_E^d$ to a perverse sheaf on $\Coh^d_0$ called the
Goresky-MacPherson extension. We denote by $\Ll_E$ the perverse sheaf
on $\Coh_0$, whose restriction to $\Coh^d_0$ equals $\Ll_E^d$.

Using the perverse sheaf $\Ll_E^d$ we define the averaging functor
$\Hav^d_{k,E}$ on the derived category of perverse sheaves on
$\Bun_k$, where the positive integer $k$ is independent of $n$, the
rank of the local system $E$. For $d\geq 0$, introduce the stack
$\Mod^d_k$, which classifies the data of triples
$(\M,\M',\beta:\M\hookrightarrow \M')$, where $\M,\M'\in \Bun_k$ and
$\beta$ is an embedding of coherent sheaves such that the quotient
$\M'/\M$ is a torsion sheaf of length $d$. We have the diagram
$$
\begin{array}{ccccc}
& &  \Mod^d_k & & \\
& \stackrel{\hl}\swarrow & & \stackrel{\hr}\searrow & \\
\Bun_k & & & & \Bun_k
\end{array}
$$
where $\hl$ (resp., $\hr$) denotes the morphism sending a triple
$(\M,\M',\beta)$ to $\M$ (resp., $\M'$). In addition, we have a
natural smooth morphism $\pi:\Mod^d_k\to \Coh^d_0$, which sends a
triple $(\M,\M',\beta)$ to the torsion sheaf $\M'/\M$.

The {\em averaging functor} $\Hav^d_{k,E}$ is defined by the formula
\begin{equation*}
\K \mapsto \hr_!(\hl{}^*(\K)\otimes \pi^*(\Ll^d_E)).
\end{equation*}
The following theorem was stated as a conjecture in \cite{FGV} and
proved in \cite{FGV} (when $X$ is defined over a finite field) and in
\cite{Ga} (when $X$ is defined over a finite field or over $\C$).

\begin{thm} \label{vanishing conjecture}
Let $E$ be an irreducible local system of rank $n$ on $X$. Then for
all $k=1,\ldots,n-1$ and all $d$ satisfying $d > kn(2g-2)$, the
functor $\Hav^d_{k,E}$ is identically equal to $0$.
\end{thm}

The geometric Langlands conjecture for $GL_n$ follows from this
theorem together with the main theorem of \cite{FGV} which states that
if \thmref{vanishing conjecture} holds for an irreducible rank $n$
local system $E$, then the geometric Langlands \conjref{glc} also
holds for $E$.

In the case when $k=1$ \thmref{vanishing conjecture} has been proved
earlier by Deligne, and this result was one of the main steps in
Drinfeld's proof of the geometric Langlands conjecture for $GL_2$
\cite{Dr1}. In this case the result may be reformulated without using
Laumon's sheaf. Recall the Abel-Jacobi map $\pi_d: S^d X \to
\on{Pic}_d$ introduced in \secref{gacft}. Then for any irreducible
rank $n>1$ local system $E$ on $X$ we have $$\pi_{d*}(E^{(d)}) = 0,
\qquad \forall d>n(2g-2).$$

It is not difficult to show that this statement is equivalent to the
following result: $H^\bullet(S^d X,(E \otimes E')^{(d)}) = 0 $ for any
rank one local system $E'$ on $X$. To prove that, observe that since
$E$ is irreducible and has rank $n>1$, so is $E \otimes E'$. Therefore
we have $H^0(X,E \otimes E') = H^2(X,E \otimes E') = 0$. But then by
the K\"unneth formula we have
$$H^\bullet(S^d X,(E \otimes E')^{(d)}) = \wedge^d H^1(X,E \otimes
E').$$ Note that $\dim H^1(X,E \otimes E') = n(2g-2)$, which
follows immediately from the computation of the Euler characteristic
of $H^\bullet(X,E_0)$, where $E_0$ is the trivial local system of rank
$n$. Therefore $$H^\bullet(S^d X,(E \otimes E')^{(d)}) = 0, \qquad
\forall d > n(2g-2),$$ which is what we needed to prove.

\section{From $GL_n$ to other reductive groups}    \label{reductive}

One adds a new dimension to the Langlands Program by considering
arbitrary reductive groups instead of the group $GL_n$. This is when
some of the most beautiful and mysterious aspects of the Program are
revealed, such as the appearance of the Langlands dual group. In this
section we will trace the appearance of the dual group in the
classical context and then talk about its
geometrization/categorification.

\subsection{The spherical Hecke algebra for an arbitrary reductive
  group}

Suppose we want to find an analogue of the Langlands correspondence
from \thmref{langl} where instead of automorphic representations of
$GL_n(\AD)$ we consider automorphic representations of $G(\AD)$, where
$G$ is a (connected, split) reductive group over $\Fq$. We wish to
relate those representations to some data corresponding to the Galois
group $G_F$. In the case of $GL_n$ this relation satisfies an
important compatibility condition that the Hecke eigenvalues of an
automorphic representation coincide with the Frobenius eigenvalues of
the corresponding Galois representation. Now we need to find an
analogue of this compatibility condition for general reductive
groups. The first step is to understand the structure of the spherical
Hecke algebra ${\mc H}_x$. For $G=GL_n$ we saw that this algebra is
isomorphic to the algebra of symmetric Laurent polynomials in $n$
variables. Now we need to give a similar description of ${\mc H}_x$
for a general group $G$.

So let $G$ be a connected reductive group over a finite field $\Fq$
which is split over $\Fq$, i.e., contains a maximal torus $T$ which is
isomorphic to a power of the multiplicative group. Then we attach to
this torus two lattices, $P$ and $P^\vee$, or characters and
cocharacters, respectively. They contain subsets $\Delta$ and
$\Delta^\vee$ of roots and coroots of $G$, respectively (see
\cite{Springer} for more details). Let us pick a point $x \in |X|$ and
assume for simplicity that its residue field is $\Fq$. To simplify
notation we will omit the index $x$ from our formulas in this
section. Thus, we will write ${\mc H}, F, \OO$ for ${\mc H}_x, F_x,
\OO_x$, etc.

The Hecke algebra ${\mc H} = {\mc H}(G(F),G(\OO))$ is by
definition the space of $\oQl$--valued compactly supported functions on
$G(F)$ which are bi-invariant with respect to the maximal compact
subgroup $G(\OO)$. It is equipped with the convolution product
\begin{equation}    \label{conv}
(f_1 \cdot f_2)(g) = \int_G f_1(x) f_2(gx^{-1}) \; dx,
\end{equation}
where $dx$ is the Haar measure on $G(F)$ normalized so that the
volume of $G(\OO)$ is equal to $1$.

We have a natural restriction homomorphism ${\mc H} \to
{\mc H}(T(F),T(\OO))$. The Hecke algebra ${\mc
H}(T(F),T(\OO))$ is easy to describe. For each $\la \in P^\vee$
we have an element $\la(t) \in T(F)$, where $t$ is a uniformizer
in $\OO$, and $T(\OO)\bs T(F)/T(\OO) = \{ \la(t) \}_{\la \in
P^\vee}$. Therefore ${\mc H}(T(F),T(\OO))$ is the group algebra
$\oQl[P^\vee]$ of $P^\vee$. The following result is called the Satake
isomorphism.

\begin{thm}    \label{satake}
The homomorphism ${\mc H} \to {\mc H}(T(F),T(\OO)) =
\oQl[P^\vee]$ is injective and its image is equal to the subalgebra
$\oQl[P^\vee]^W$ of $W$--invariants, where $W$ is the Weyl group of $G$.
\end{thm}

A crucial observation of R. Langlands was that $\oQl[P^\vee]^W$ is
nothing but the representation ring of the group $^L G(\oQl)$, the
{\em Langlands dual group} of $G$ \cite{Lan}. By definition, $^L
G(\oQl)$ is the reductive group over $\oQl$ with a maximal torus $^L
T(\oQl)$ dual to $T$, so that its lattices of characters and
cocharacters are those of $T$ interchanged, and with the sets of roots
and coroots being those of $G$, also interchanged (see
\cite{Springer}). For instance, the dual of $GL_n$ is again $GL_n$,
$SL_n$ is dual to $PGL_n$, $SO_{2n+1}$ is dual to $Sp_n$, and
$SO_{2n}$ is self-dual.

Let $\Rep {}^L G$ be the Grothendieck ring of the category of
finite-dimensional representations of $^L G(\oQl)$. Then the character
homomorphism $\Rep {}^L G \to \oQl[P^\vee]$ is injective and its image is
equal to $\oQl[P^\vee]^W$. Therefore \thmref{satake} may be interpreted
as saying that ${\mc H} \simeq \Rep {}^L G$. It follows then that the
one-dimensional representations of ${\mc H}$ are nothing but the
semi-simple conjugacy classes of $^L G(\oQl)$. Indeed, if $\ga$ is a
semi-simple conjugacy class in $^L G(\oQl)$, then we attach to it a
one-dimensional representation of ${\mc H}$ by the formula $[V]
\mapsto \on{Tr}(\ga,V)$.

\subsection{Towards the Langlands correspondence for an arbitrary
reductive group}

Now we can formulate for an arbitrary reductive group $G$ an analogue
of the compatibility statement in the Langlands conjecture
\conjref{langl} for $GL_n$. Namely, suppose that $\pi = \bigotimes'_{x
\in |X|} \pi_x$ is a cuspidal automorphic representation of
$G(\AD)$. For all but finitely many $x \in |X|$ the representation
$\pi_x$ of $G(F_x)$ is unramified, i.e., the space of
$G(\OO_x)$--invariants in $\pi_x$ is non-zero. Then this space of
$G(\OO_x)$--invariants is one-dimensional and so ${\mc H}_x$ acts on
it via a character, which by \thmref{satake} corresponds to a
semi-simple conjugacy class $\ga_x$ in $^L G(\oQl)$. Thus, we attach
to an automorphic representation a collection $\{ \ga_x \}$ of
semi-simple conjugacy classes in $^L G(\oQl)$ for almost all points of
$X$. Therefore on the other side of the Langlands correspondence we
need some sort of Galois data which would also give us such a
collection of conjugacy classes.

The candidate that immediately comes to mind is a homomorphism
$$\sigma: G_F \to {}^L G(\oQl),$$ which is almost everywhere
unramified. Then we may attach to $\sigma$ a collection of conjugacy
classes $\{ \sigma(\on{Fr}_x) \}$ of $^L G(\oQl)$ at almost all points
of $X$. So in the first approximation we may formulate the Langlands
correspondence for general reductive groups as a correspondence
between automorphic representations of $G(\AD)$ and Galois
homomorphisms $G_F \to {}^L G(\oQl)$ which satisfies the following
compatibility condition: if $\pi$ corresponds to $\sigma$, then the
$^LG$--conjugacy classes attached to $\pi$ through the action of the
Hecke algebra are the same as the Frobenius $^LG$--conjugacy classes
attached to $\sigma$.

Unfortunately, the situation is not as clear-cut as in the case of
$GL_n$ because many of the results which facilitate the Langlands
correspondence for $GL_n$ are no longer true in general. For instance,
it is not true that the collection of the Hecke conjugacy classes
determines the automorphic representation uniquely or that the
collection of the Frobenius conjugacy classes determines the Galois
representation uniquely in general. Therefore even the statement of
the Langlands conjecture becomes a much more subtle issue for a
general reductive group. However, the main idea appears to be correct:
there is a relationship, still very mysterious, alas, between
automorphic representations of $G(\AD)$ and homomorphisms from the
Galois group $G_F$ to the Langlands dual group $^L G$.

Now (in the hope of gaining some insight into this mystery) we would
like to formulate a geometric analogue of this relationship. The first
step is to develop a geometric version of the Satake isomorphism.

\subsection{Categorification of the spherical Hecke algebra}
\label{cate}

Let us look at the isomorphism of \thmref{satake} more closely. It is
easy to see that the elements $\la(t)$, where $\la \in P^\vee_+$, the
set of dominant weights (with respect to a Borel subgroup of $^L G$),
are representatives of the double cosets of $G(F)$ with respect to
$G(\OO)$. Therefore ${\mc H}$ has a basis $\{ c_\la \}_{\la \in
P^\vee_+}$, where $c_\la$ is the characteristic function of the double
coset $G(\OO) \la(t) G(\OO) \subset G$. On the other hand, $\Rep {}^L
G$ also has a basis labeled by $\la \in P^\vee_+$, which consists of
the classes $[V_\la]$, where $V_\la$ is the irreducible representation
with highest weight $\la$. However, under the Satake isomorphism these
bases do not coincide. Instead, we have the following formula
\begin{equation}    \label{Hla}
H_\la =  q^{-(\la,\rho)} \left( c_\la + \sum_{\mu \in P^\vee_+;
\mu<\la} a_{\la\mu} c_\mu \right), \quad \quad a_{\la\mu} \in \Z_+[q],
\end{equation}
where $H_\la$ is the image of $[V_\la]$ in ${\mc H}$ under the
Satake isomorphism. This formula has a remarkable geometric
explanation.

Let us consider ${\mc H}$ as the algebra of functions on the
quotient $G(F)/G(\OO)$ which are left invariant with respect to
$G(\OO)$. In view of the Grothendieck fonctions-faisceaux dictionary
discussed in \secref{grot}, it is natural to ask whether
$G(F)/G(\OO)$ is the set of $\Fq$--points of an algebraic variety, and
if so, whether $H_\la$ is the function corresponding to a perverse
sheaf on this variety. It turns out that this is indeed the case.

The quotient $G(F)/G(\OO)$ is the set of points of an ind-scheme $\Gr$
over $\Fq$ called the {\em affine Grassmannian} associated to $G$. Let
${\mc P}_{G(\OO)}$ be the category of $G(\OO)$--equivariant (pure)
perverse sheaves on $\on{Gr}$. For each $\la \in P^\vee_+$ we have a
finite-dimensional $G(\OO)$--orbit $\Gr_\la = G(\OO) \cdot \la(t)$ in
$\Gr$. Let $\IC_\la$ be the irreducible perverse sheaf obtained by the
Goresky-MacPherson extension from the constant sheaf on $\Gr_\la$ to
its closure $\ol\Gr_\la$. These are the irreducible objects of the
category ${\mc P}_{G(\OO)}$.

Assigning to a perverse sheaf its ``trace of Frobenius'' function, we
obtain an identification between the Grothendieck group of ${\mc
P}_{G(\OO)}$ and the algebra of $G(\OO)$--invariant functions on
$G(F)/G(\OO)$, i.e., the spherical Hecke algebra. In that sense, ${\mc
P}_{G(\OO)}$ is a {\em categorification} of the Hecke algebra. A
remarkable fact is that the function $H_\la$ in formula \eqref{Hla} is
equal to the function associated to the perverse sheaf $\IC_\la$, up
to a sign $(-1)^{2(\la,\rho)}$. Thus, under the Satake isomorphism the
classes of irreducible representations of $^L G$ go not to functions
which correspond to constant sheaves on the orbits (i.e., the
functions $c_\la$) but to the irreducible {\em perverse} sheaves. This
suggests that the Satake isomorphism itself may be elevated from the
level of Grothendieck groups to the level of categories. This is
indeed true.

In fact, it is possible to define the structure of tensor category on
${\mc P}_{G(\OO)}$ with the tensor product given by a convolution
functor corresponding to the convolution product \eqref{conv} at the
level of functions. Then up to a small subtlety, which has to do with
the appearance of the sign $(-1)^{2(\la,\rho)}$ mentioned above, we
have the following beautiful result. It has been conjectured by
Drinfeld and proved by I. Mirkovi\'c and K. Vilonen \cite{MV} and
V. Ginzburg \cite{Ginzburg} (some important results in this direction
were obtained earlier by G. Lusztig \cite{Lusztig}).

\begin{thm}
The tensor category ${\mc P}_{G(\OO)}$ is equivalent to the tensor
category of finite-dimensional representations of the group $^L
G(\oQl)$.
\end{thm}

Moreover, the fiber functor from ${\mc P}_{G(\OO)}$ to the category of
vector spaces (corresponding to the forgetful functor from the
category of representations) is the global cohomology functor ${\mc F}
\mapsto \oplus_i H^i(\Gr,{\mc F})$. This allows one to reconstruct the
Langlands dual group $^L G$ by means of the standard Tannakian
formalism. So we get a completely new perspective on the nature of the
dual group. For example, the dual group to $GL_n$ now appears as the
group of automorphisms of the total cohomology space of the
projectivization of its $n$--dimensional defining representation.
This is a good illustration of why geometry is useful in the Langlands
Program.

The above theorem should be viewed as a categorification of the Satake
isomorphism of \thmref{satake}. We can use it to define the notion of
a Hecke eigensheaf for an arbitrary reductive group and to formulate a
geometric version of the Langlands correspondence. In the next section
we do that for curves over $\C$, but one can apply the same technique
over the finite field as well.

\section{The geometric Langlands conjecture over $\C$}    \label{over
C}

In this section we will formulate the geometric Langlands conjecture
for an arbitrary reductive group over $\C$. We will then give a brief
overview of the recent work of A. Beilinson and V. Drinfeld in which a
substantial part of this conjecture has been proved. It is interesting
that their work uses results from representation theory of affine
Kac-Moody algebras \cite{FF:gd,F:wak}, which now play the role of the
reductive groups over local non-archimedian fields.

\subsection{Hecke eigensheaves}

In the rest of this paper $X$ will be a smooth connected projective
curve over $\C$ and $G$ a reductive algebraic group over $\C$. The
results of \secref{cate} are applicable in this context. Namely, we
have the affine Grassmannian over $\C$ and the category ${\mc
P}_{G(\OO)}$ of $G(\OO)$--equivariant perverse sheaves (of
$\C$--vector spaces) on $\Gr$. This category is equivalent, as a
tensor category, to the category of finite-dimensional representations
of the Langlands dual group $^L G(\C)$. Under this equivalence, the
irreducible representation of $^L G$ with highest weight $\la \in
P^\vee_+$ corresponds to the irreducible perverse sheaf $\IC_\la$.

Now we can define the analogues of the $GL_n$ Hecke functors
introduced in \secref{Hecke functors} for general reductive
groups. Let $\Bun_G$ be the moduli stack of $G$--bundles on
$X$. Consider the stack ${\mc H}ecke$ which classifies quadruples
$(\M,\M',x,\beta)$, where $\M$ and $\M'$ are $G$--bundles on $X$, $x
\in X$, and $\beta$ is an isomorphism between the restrictions of $\M$
and $\M'$ to $X \bs x$. We have natural morphisms
$$
\begin{array}{ccccc}
& & {\mc Hecke} & & \\
& \stackrel{\hl}\swarrow & & \stackrel{\hr}\searrow & \\
\Bun_G & & & & X\times \Bun_G
\end{array}
$$
where $\hl(\M,\M',x,\beta) = \M$ and $\hr(\M,\M',x,\beta) = (x,\M')$.

Note that the fiber of ${\mc H}ecke$ over $(x,\M')$ is the moduli
space of pairs $(\M,\beta)$, where $\M$ is a $G$--bundles on $X$, and
$\beta: \M'|_{X\bs x} \overset{\sim}\to \M|_{X\bs x}$. It is known
that this is a twist of $\Gr_x = G(F_x)/G(\OO_x)$ by the
$G(\OO)_x$--torsor $\M'(\OO_x)$ of sections of $\M'$ over $\on{Spec}
\OO_x$:
$$
(\hr)^{-1}(x,\M') = \M'(\OO_x) \underset{G(\OO_x)}\times \Gr_x.
$$
Therefore we have a stratification of each fiber, and hence of the
entire ${\mc H}ecke$, by the substacks ${\mc H}ecke_\la, \la \in
P^\vee_+$, which correspond to the $G(\OO)$--orbits $\Gr_\la$ in
$\Gr$. Consider the irreducible perverse sheaf on ${\mc H}ecke$, which
is the Goresky-MacPherson extension of the constant sheaf on ${\mc
H}ecke_\la$. Its restriction to each fiber is isomorphic to $\IC_\la$,
and by abuse of notation we will denote this entire sheaf also by
$\IC_\la$.

Define the Hecke functor $\He_\la$ from the derived category of
perverse sheaves on $\Bun_G$ to the derived category of perverse
sheaves on $X \times \Bun_G$ by the formula
$$
\He_\la({\mc F}) = \hr_!(\hl{}^*({\mc F}) \otimes \IC_\la).
$$
Let $E$ be a $^L G$--local system on $X$. Then for each irreducible
representation $V_\la$ of $^L G$ we have a local system $V_\la^E = E
\underset{G}\times V_\la$.

\begin{definition}
{\em A perverse sheaf on $\Bun_G$ is a called a {\em Hecke eigensheaf
with eigenvalue} $E$ if we are given isomorphisms
$$
\He_\la({\mc F}) \simeq V_\la^E \boxtimes {\mc F},
$$
which are compatible with the tensor product structure on the category
of representations of $^L G$.}
\end{definition}

In the case when $G=GL_n$ this definition is equivalent to equations
\eqref{eigen-property}.

Now we can state the geometric Langlands conjecture.

\begin{conj}    \label{glc1}
Let $E$ be a $^L G$--local system on $X$ which cannot be reduced to a
proper parabolic subgroup of $^L G$. Then there exists a non-zero
Hecke eigensheaf $\Aut_E$ on $\Bun_G$ with the eigenvalues $E$ which is
irreducible and perverse on each connected component of $\Bun_G$.
\end{conj}

When working over $\C$, we may switch from perverse sheaves to
$\D$--{\em modules}. If $V$ is a smooth variety, we consider the sheaf
of differential operators on $V$ and sheaves of modules over it, which
we simply refer to as $\D$--modules. The simplest example of a
$\D$--module is the sheaf of sections of a vector bundle on $V$
equipped with a flat connection (we can use the flat connection to act
on sections by vector fields). The sheaf of horizontal sections of
this bundle is then a locally constant sheaf, and hence a perverse
sheaf. One can associate to a more general $\D$--module a perverse
sheaf in a similar way. In fact, there is an equivalence between the
category of holonomic $\D$--modules with regular singularities on a
variety $V$ and the category of perverse sheaves on $V$, called the
Riemann-Hilbert correspondence (see \cite{Dmodules,GM}). Therefore we
may replace in the above conjecture perverse sheaves by
$\D$--modules. In what follows we will consider this $\D$--module
version of the geometric Langlands conjecture.

\subsection{The geometric Langlands correspondence as a Fourier-Mukai
transform}    \label{fm}

Consider first the case of $G=GL_1$. Then $\Bun_1$ is the Picard
variety $\on{Pic}$ and in order to prove the conjecture we should
attach a Hecke eigensheaf $\Aut_E$ on $\on{Pic}$ to each rank one
local system $E$ on $X$. Let $\sigma: \pi_1(X) \to \C^\times$ be the
homomorphism corresponding to $\sigma$. It factors through the maximal
abelian quotient of $\pi_1(X)$, i.e., $H_1(X,\Z)$. But using the
cup-product on $H_1(X,\Z)$ we may identify it with $H^1(X,\Z)$, which
is equal to $\pi_1(\on{Pic}_0)$. Therefore $\sigma$ gives rise to a
homomorphism $\pi_1(\on{Pic}_0) \to \C^\times$ and hence to a rank one
local system on $\on{Pic}_0$, or equivalently, a flat line bundle on
$\on{Pic}_0$. This flat line bundle (considered as a $\D$--module) is
precisely the restriction of $\Aut_E$ to $\on{Pic}_0$. It may be
extended to the other components of $\on{Pic}$ in the same way as in
\secref{gacft}.

So in the case of $GL_1$ the Hecke eigensheaves are actually flat line
bundles. One constructs them by using the fact that rank one local
systems on a curve $X$ are the same as rank one local systems on its
Jacobean $\on{Jac} = \on{Pic}_0$. Note however that we have used the
isomorphism $H_1(X,\Z) \simeq H^1(X,\Z)$. Because of that, if we take
$G$ to be an arbitrary torus $T$, then rank one local systems on
connected components of $\Bun_T$ will correspond to $^L T$--local
systems on $X$.

One can strengthen the statement of the geometric Langlands conjecture
by using the Fourier-Mukai transform. Let $\on{Loc}_1$ be the moduli
space of rank one local systems on $X$, or equivalently, on
$\on{Jac}$. On the product $\on{Loc}_1 \times \on{Jac}$ we have the
``universal flat line bundle'' $\Aut$, whose restriction to $\{ E \}
\times \on{Jac}$ is the flat line bundle on $\on{Jac}$ corresponding
to $E$. It has a partial flat connection along $\on{Jac}$. This
enables us to define functors between the derived category of
$\OO$--modules on $\on{Loc}_1$ and the derived category of ${\mc
D}$--modules on $\Bun_1$: pulling back to $\on{Loc}_1 \times
\on{Jac}$, tensoring with $\Aut$ and pushing forward to the other
factor. It has been shown by G. Laumon \cite{Laumon:F} and
M. Rothstein \cite{Rothstein} that these functors give rise to
mutually inverse equivalences of derived categories, up to a sign. We
note that a similar equivalence holds if one replaces $\on{Jac}$ by an
arbitrary abelian variety $A$, and it generalizes the original
Fourier-Mukai correspondence in which one considers $\OO$--modules on
$A$ as opposed to $\D$--modules.

Note that under this equivalence the skyscraper sheaf supported at the
point $\{ E \} \in \on{Loc}_1$ goes to the sheaf $\Aut_E$ on
$\on{Jac}$. So the above equivalence may be loosely interpreted as
saying that any $\D$--module on $\on{Jac}$ may be expressed as a
``direct integral'' of the Hecke eigensheaves $\Aut_E$. In other
words, the Fourier-Mukai equivalence may be viewed as a ``spectral
decomposition'' of the derived category of $\D$-modules on $\on{Jac}$.

Optimistically, one may hope that the Langlands correspondence for
general reductive groups also gives a kind of spectral decomposition
of the derived category of $\D$-modules on $\Bun_G$ (or, more
precisely, its connected component). Namely, one may hope that there
exists an equivalence between this derived category and the derived
category of $\OO$--modules on the moduli stack $\on{Loc}_{^L G}$ of
$^L G$--local systems on $X$, so that the skyscraper sheaf on
$\on{Loc}_{^L G}$ supported at the local system $E$ corresponds to the
Hecke eigensheaf $\Aut_E$. If this is true, it would mean that Hecke
eigensheaves provide a good ``basis'' in the category of $\D$--modules
on $\Bun_G$, just as flat line bundles provide a good ``basis'' in the
category of $\D$--modules on $\on{Jac}$.

While it is not known whether this non-abelian Fourier-Mukai transform
exists, A. Beilinson and V. Drinfeld have recently constructed a
special case of this transform for an arbitrary semisimple group
$G$. Roughly speaking, they construct $\D$--modules on $\Bun_G$
corresponding to $\OO$--modules supported on a certain affine
subvariety in $\on{Loc}_{^L G}$ called the space of $^L
G$--opers. Before discussing their construction, we consider its
analogue in the abelian case.

\subsection{A special case of Fourier-Mukai transform}
\label{special}

The moduli space $\on{Loc}_1$ of flat line bundles on $X$ fibers over
$\on{Jac} = \on{Pic}_0$ with the fiber over $\Ll \in \on{Jac}$ being
the space of all (holomorphic) connections on $\Ll$. This is an affine
space over the space $H^0(X,\Omega)$ of holomorphic one-forms on
$X$. In particular, the fiber over the trivial line bundle $\OO$ is
just the space $H^0(X,\Omega)$. As we have seen above, each point
$\omega \in H^0(X,\Omega)$ gives rise to a flat line bundle on
$\on{Jac}$. It turns out that we can describe the $\D$--module of
sections of this flat line bundle quite explicitly.

Observe that because the tangent bundle to $\on{Jac}$ is trivial, with
the fiber isomorphic to $H^1(X,\Omega)$, the algebra $D$ of global
differential operators on $\on{Jac}$ is commutative and is isomorphic
to $\on{Sym} H^1(X,\OO) = \on{Fun} H^0(X,\Omega)$, by the Serre
duality. Therefore each point $\omega \in H^0(X,\Omega)$ gives rise to
a character $\chi_\omega: D \to \C$. Define the $\D$--module ${\mc
F}_\omega$ on $\on{Jac}$ by the formula
\begin{equation}    \label{abelian}
{\mc F}_\omega = {\mc D}/I_\omega,
\end{equation}
where ${\mc D}$ is the sheaf of differential operators on $\on{Jac}$,
and $I_\omega$ is the left ideal in ${\mc D}$ generated by the kernel
of $\chi_\omega$ in $D$. Then this is a holonomic $\D$--module which
is equal to the Hecke eigensheaf corresponding to the trivial line
bundle with connection $d+\omega$ on $X$.

We note that the $\D$--module ${\mc F}_\omega$ represents the system
of differential equations
\begin{equation}    \label{eqs}
X \cdot \Psi = \chi_\omega(X) \Psi, \qquad X \in D,
\end{equation}
in the sense that for any homomorphism from ${\mc F}_\omega$ to
another $\D$--module $\K$ the image of $1 \in {\mc F}_\omega$ in $\K$
is a solution of the system \eqref{eqs}.

Generalizing the definition of the $\D$--module ${\mc F}_\omega$ we
obtain a functor from the category of modules over $\on{Fun}
H^0(X,\Omega)$ to the category of $\D$--modules on $\on{Jac}$,
\begin{equation}    \label{system}
M \mapsto \D \underset{D}\otimes M,
\end{equation}
so that $\chi_\omega \mapsto {\mc F}_\omega$. This functor is the
restriction of the Fourier-Mukai functor to the category of
$\OO$--modules on $\on{Loc}_1$ supported on $H^0(X,\Omega) \subset
\on{Loc}_1$.

Beilinson and Drinfeld have given a similar construction for an
arbitrary semisimple group $G$.

\subsection{Opers}

Let $G$ be a simple algebraic group, which we will assume to be
connected and simply-connected in the rest of this section. Then the
dual group $^L G$ is of adjoint type. The analogue of the affine space
of connections on the trivial line bundle is the space of $^L
G$--opers.

Let $G$ be a simple algebraic group of adjoint type, $B$ a Borel
subgroup and $N = [B,B]$ its unipotent radical, with the corresponding
Lie algebras $\n \subset \bb\subset \g$. There is an open $B$--orbit
${\bf O}\subset \n^\perp/\bb \subset \g/\bb$, consisting of vectors
which are stabilized by the radical $N\subset B$, and such that all of
their negative simple root components, with respect to the adjoint
action of $H = B/N$, are non-zero. This orbit may also be described as
the $B$--orbit of the sum of the projections of simple root generators
$f_i$ of any nilpotent subalgebra $\n_-$, which is in generic position
with $\bb$, onto $\g/\bb$. The torus $H = B/N$ acts simply
transitively on ${\bf O}$, so ${\bf O}$ is an $H$--torsor.

Suppose we are given a principal $G$--bundle $\F$ on a smooth curve
$X$ with a connection $\nabla$ and a reduction $\F_B$ to the Borel
subgroup $B$ of $G$. Then we define the relative position of $\nabla$
and $\F_B$ (i.e., the failure of $\nabla$ to preserve $\F_B$) as
follows. Locally, choose any flat connection $\nabla'$ on $\F$
preserving $\F_B$, and take the difference $\nabla - \nabla'$.  It is
easy to show that the resulting local sections of $(\g/\bb)_{\F_B}
\otimes \Omega$ are independent of $\nabla'$, and define a global
$(\g/\bb)_{\F_B}$--valued one-form on $X$, denoted by $\nabla/\F_B$.

\begin{definition}
{\em A $G$--{\em oper} on $X$ is by definition a triple
$(\F,\nabla,\F_B)$, where $\F$ is a principal $G$--bundle $\F$ on $X$,
$\nabla$ is a connection on $\F$ and $\F_B$ is a $B$--reduction of
$\F$, such that the one--form $\nabla/\F_B$ takes values in ${\bf
O}_{\F_B} \subset(\g/\bb)_{\F_B}$.}
\end{definition}

This definition is due to Beilinson and Drinfeld \cite{BD} (in the
case when $X$ is replaced by the punctured disc $\on{Spec} \C((t))$
opers were introduced earlier by Drinfeld and Sokolov in their work on
the generalized KdV hierarchies \cite{DS}). Note that ${\bf O}$ is
$\C^\times$--invariant, so that ${\bf O}_{\F_B}$ is well-defined in
$(\g/\bb)_{\F_B}$.

For instance, in the case when $G=PGL_n$ this condition means that if
we choose a local trivialization of $\F_B$ and a local coordinate $t$
on $X$, then the connection operator will have the form
$$
\nabla = \pa_t + \left( \begin{array}{ccccc}
*&*&*&\hdots&*\\
\star&*&*&\hdots&*\\
0&\star&*&\hdots&*\\
\vdots&\ddots&\ddots&\ddots&\vdots\\ 
0&0&\hdots&\star&*
\end{array} \right)
$$
where the $*$'s indicate arbitrary functions in $t$ and the $\star$'s
indicate nowhere vanishing functions.

By changing the trivialization of $\F_B$ this operator may be brought
in a unique way to the form
$$
\partial_t + \left( \begin{array}{ccccc}
0&v_1&v_2&\cdots&v_{n-1}\\
-1&0&0&\cdots&0\\
0&-1&0&\cdots&0\\
\vdots&\ddots&\ddots&\cdots&\vdots\\
0&0&\cdots&-1&0
\end{array}\right).
$$
But giving such an operator is the same as giving a scalar $n$th
order scalar differential operator
\begin{equation}    \label{first time opers}
\partial_t^n+v_1(t) \partial_t^{n-2}+\ldots+v_{n-1}(t),
\end{equation}
and a more careful calculation shows that it must act from
$\Omega^{-(n-1)/2}$ to $\Omega^{(n+1)/2}$). So the space of
$PGL_n$--opers is the space of operators of the form \eqref{first time
opers} (if $n$ is even, we need to choose a square root of the
canonical line bundle $\Omega$, but the space of $PGL_n$--opers is
independent of this choice). In particular, it turns out that a
$PGL_2$--oper is nothing but a projective connection.

\begin{lem}
If $G$ is a simple algebraic group of adjoint type and $X$ is a smooth
projective curve of genus not equal to one, then there exists a unique
(up to isomorphism) $G$--bundle $\F_0$ which admits the structure of
an oper. Moreover, the corresponding Borel reduction $\F_{0,B}$ is
also uniquely determined, and for any connection $\nabla$ on $\F_0$
the triple $(\F_0,\nabla,\F_{0,B})$ is a $G$--oper.
\end{lem}

For example, for $G=PGL_n$ the bundle $\F_0$ corresponds to the rank
$n$ vector bundle of $(n-1)$--jets of sections of the line bundle
$\Omega^{-(n-1)/2}$ (note that the corresponding $PGL_n$--bundle is
independent of the choice of $\Omega^{1/2}$).

Now we switch to the Langlands dual group. The above lemma means that
the space $\on{Op}_{^L G}(X)$ of $^L G$--opers on $X$ is an affine
space which is identified with the fiber of the natural projection
$\on{Loc}_{^L G} \to \Bun_{^L G}$ over $\F_0$.

Beilinson and Drinfeld associate to each point of $\on{Op}_{^L G}(X)$
a Hecke eigensheaf on $\Bun_G$. Recall that in the abelian case the
crucial point was that the algebra of global differential operators on
$\Jac$ was isomorphic to $\on{Fun} H^0(X,\Omega)$. Beilinson and
Drinfeld prove an analogue of this statement in the non-abelian
case. However, in this case it is necessary to consider the sheaf
$\D'$ of differential operators acting on the square root of the
canonical line bundle on $\Bun_G$ (this square root is unique up to
isomorphism).

\begin{thm}    \label{bd}
The algebra of global sections of the sheaf $\D'$ is commutative and
is isomorphic to the algebra of functions on the space $\on{Op}_{^L
G}(X)$ of $^L G$--opers on $X$.
\end{thm}

\subsection{Hitchin's integrable system}

It is instructive to look at the quasi-classical analogue of this
statement. The algebra $D' = H^0(\Bun_G,\D')$ carries the standard
filtration by the order of the differential operator, and the
associated graded algebra embeds into the algebra of functions on the
cotangent bundle $T^* \Bun_G$ to $\Bun_G$. On the other hand, it is
not difficult to show that $\on{Op}_{^L G}(X)$ is an affine space over
the space
$$
H_G(X) = \bigoplus_{i=1}^\ell H^0(X,\Omega^{\otimes(d_i+1)}),
$$
where $\ell$ is the rank of $G$, and the $d_i$'s are the exponents of
$G$. Therefore the algebra $\on{Fun} \on{Op}_{^L G}(X)$ carries a
filtration such that the associated graded is $\on{Fun} H_G$. The
quasi-classical analogue of the isomorphism of \thmref{bd} is an
isomorphism $\on{Fun} T^* \Bun_G \simeq \on{Fun} H_G$.

To construct such an isomorphism, we need a morphism $p: T^* \Bun_G
\to H_G$. Such a morphism has been constructed by N. Hitchin
\cite{Hitchin}. Namely, let $\g$ be the Lie algebra of $G$. It is
well-known that the algebra of invariant functions on $\g^*$ is
isomorphic to the graded polynomial algebra $\C[P_1,\ldots,P_\ell]$,
where $\deg P_i = d_i+1$. Let us observe that the tangent space to
$\Bun_G$ at $\M \in \Bun_G$ is isomorphic to $H^1(X,\g_{\M})$, where
$\g_{\M} = \M \underset{G}\times \g$. Hence the cotangent space at
$\M$ is isomorphic to $H^0(X,\g^*_{\M} \otimes \Omega)$ by the Serre
duality.

By definition, the Hitchin map $p$ takes $(\M,\eta) \in T^* \Bun_G$,
where $\eta \in H^0(X,\g^*_{\M} \otimes \Omega)$ to
$(P_1(\eta),\ldots,P_\ell(\eta)) \in H_G$. It has been proved in
\cite{Hitchin,Faltings} that over an open dense subset of $H_G$ the
morphism $p$ is smooth and its fibers are proper. Therefore we obtain
an isomorphism $\on{Fun} T^* \Bun_G \simeq \on{Fun} H_G$. Moreover,
for any $\phi, \psi \in \on{Fun} H_G$, we have $\{ p^*\phi,p^*\psi \}
= 0$, where $\{\cdot,\cdot\}$ is the natural Poisson structure on $T^*
\Bun_G$ (so that $p$ gives rise to an algebraic completely integrable
system). This is a precursor of the commutativity property of the
global differential operators.

\subsection{Beilinson--Drinfeld construction}

How can we ``quantize'' the map $p$, i.e., construct an algebra
homomorphism $\on{Fun} \on{Op}_{^L G}(X) \to D'$?

In order to do this Beilinson and Drinfeld apply the following general
construction. Suppose $\wt{M} = \on{Spec} \wt{B}$ is an affine
algebraic variety with an action of an algebraic group $H$. Let $K$ be
a connected subgroup of $H$ and $M = \on{Spec} B$ another affine
variety such that $M=\wt{M}/K$, so that $B = \wt{B}^K =
\wt{B}^{\mathfrak k}$, where $\kk = \on{Lie} K$. Set $I = U\h \cdot
\kk$, where $\h = \on{Lie} H$, and $$N(I) = \{ a \in U\h \; | \; Ia
\subset I \}.$$ Then $N(I)/I$ acts on $B$ and hence we obtain a
homomorphism $N(I)/I \to H^0(M,\D_M)$.

The algebra $N(I)/I$ may also be described as follows. Consider the
induced $\h$--module $V = \on{Ind}_{\kk}^{\h} \C$. Then $$N(I)/I =
V^{\kk} = (\on{End}_{\h} V)^{\on{opp}}.$$

In our case we let $M$ be $\Bun_G$. Let $x$ be a point of $X$. In the
case when the group $G$ is semisimple, any $G$--bundle on $X$ may be
trivialized already on $X\bs x$ and so one has the following stronger
version of \lemref{weil}:
$$
\Bun_G \simeq G(\C[X\bs x])\bs G(F_x)/G(\OO_x).
$$
Let $$\wt{\Bun}_G = G(\C[X\bs x])\bs G(F_x)$$ be the moduli space of
pairs $(\M,s)$, where $\M$ is a $G$--bundle on $X$ and $s$ is its
trivialization on the formal disc $D_x = \on{Spec} \OO_x$. We may then
set $\wt{M} = \wt{\Bun}_G$, $H = G(F_x)$ and $K=G(\OO_x)$.  Though
$\Bun_G$ is an algebraic stack, the above general construction is
still applicable (see \cite{BD}), and so we obtain a homomorphism
$$
(\on{End}_{\g(F_x)} V)^{\on{opp}} \to H^0(\Bun_G,\D).
$$

Unfortunately, $\on{End}_{\g(F_x)} V = \C$, so we cannot obtain any
non-trivial global differential operators on $\Bun_G$. But the Lie
algebra $\g(F_x) \simeq \g((t))$ has a one-dimensional universal
central extension $\ghat_\ka$ called the {\em affine Kac-Moody
algebra}:
$$
0 \to \C K \to \ghat_\ka \to \g(F_x) \to 0
$$
(see \cite{Kac}). As a vector space it splits into a direct sum
$\g(F_x) \oplus \C K$, and the commutation relations read
$$
[A \otimes f,B \otimes g] = [A,B] \otimes f g - (\kappa(A,B)
\on{Res} f dg) K, \qquad [K,\cdot]=0.
$$
where $\ka$ is a non-degenerate invariant inner product on $\g$ (which
is unique up to a scalar).

Consider the $\ghat_\ka$--module $V_\ka = \on{Ind}_{\g(\OO_x) \oplus
\C K}^{\ghat_\ka} \C_1$, where $\C_1$ is the one-dimensional
representation of $\g(\OO_x) \oplus \C K$ on which $\g(\OO_x)$ acts by
$0$ and $K$ acts as the identity. Then one can generalize the above
construction and obtain a homomorphism from
$$
(\on{End}_{\ghat_\ka}
V_\ka)^{\on{opp}} \to \Gamma(\Bun_G,\D_\ka),
$$
where $\D_\ka$ is the sheaf of twisted differential operators
corresponding to $\ka$. In particular, if we choose the inner product
$\ka_c$ defined by the formula
$$
\ka_c(A,B) = - \frac{1}{2} \on{Tr}_\g \on{ad} A \on{ad} B,
$$
then $\D_{\ka_c} = \D'$, the sheaf of differential operators acting on
the square root of the canonical line bundle on $\Bun_G$.

The following result is due to B. Feigin and the author
\cite{FF:gd,F:wak}.

\begin{thm}
If $\ka \neq \ka_c$, then $\on{End}_{\ghat_{\ka}} V_{\ka} =\C$. If
$\ka=\ka_c$, then there is a canonical algebra isomorphism
$$
\on{End}_{\ghat_{\ka_c}} V_{\ka_c} \simeq \on{Fun} \on{Op}_{^L
G}(D_x),
$$
where $D_x = \on{Spec} \OO_x$.
\end{thm}

Thus, we obtain a homomorphism $\varphi_x: \on{Fun} \on{Op}_{^L
G}(D_x) \to D'$, where $D'$ is the algebra of global differential
operators acting on the square root $\omega^{1/2}$ of the canonical
line bundle on $\Bun_G$. Beilinson and Drinfeld prove the following
theorem \cite{BD}.

\begin{thm}
The homomorphism $\varphi_x$ factors through a homomorphism $$\varphi:
\on{Fun} \on{Op}_{^L G}(X) \to D'$$ which is independent of $x$ and is
an algebra isomorphism.
\end{thm}

This proves \thmref{bd}. Now, given a $^L G$--oper $\rho =
(\F,\nabla,\F_B)$ on $X$, we construct a $\D$--module $\Delta_\rho$ on
$\Bun_G$:
$$
\Delta_\rho = (\D'/I_\rho) \otimes \omega^{-1/2},
$$ where $I_\rho$ is the left ideal in $\D'$ generated by the kernel
of the character $\chi_\rho: D' \to \C$ corresponding to the point
$\rho \in \on{Spec} D'$ (compare with formula \eqref{abelian}). Since
$\dim \on{Op}_{^L G}(X)$ $= \dim \Bun_G$, the $\D$--module
$\Delta_\rho$ is holonomic (it is also non-zero). Denote by $E_\rho$
the $^L G$--local system $(\F,\nabla)$ underlying the $^L G$--oper
$\rho$. Beilinson and Drinfeld prove the following fundamental result.

\begin{thm}
The $\D$--module $\Delta_\rho$ is a Hecke eigensheaf with eigenvalue
$E_\rho$.
\end{thm}

Thus, Beilinson and Drinfeld prove the geometric Langlands Conjecture
for those $^L G$--local systems on $X$ which admit the structure of an
oper.

More generally, as in \secref{special} we obtain a functor from the
category of modules over the algebra $\on{Fun} \on{Op}_{^L G}(X)$ to
the category of $\D$--modules on $\Bun_G$:
$$
M \mapsto (\D' \underset{D'}\otimes M) \otimes \omega^{-1/2},
$$ so that $\chi_\rho \mapsto \Delta_\rho$. This functor may be viewed
as the restriction of the would-be non-abelian Fourier-Mukai transform
to the category of $\OO$--modules on $\on{Loc}_{^L G}$ supported on
$\on{Op}_{^L G}(X) \subset \on{Loc}_{^L G}$.

\end{document}